\documentclass[11pt, reqno]{amsart}

\usepackage{amsmath, amssymb, comment, amsthm, amssymb}
\usepackage{amscd}
\usepackage{graphicx}
\usepackage{ragged2e}
\usepackage{epstopdf}
\usepackage{setspace}
\usepackage{abstract}
\usepackage{subfigure}
\usepackage{color}
\usepackage[foot]{amsaddr}
\usepackage{enumitem}
\usepackage{tabularx}
\usepackage[percent]{overpic}

\usepackage{tikz}
\usetikzlibrary{decorations.markings}

\usepackage[format=plain,labelfont=bf,up,width=.9\textwidth]{caption}

\usepackage[linktocpage=true]{hyperref}

\theoremstyle{plain}
\newtheorem{thm}{Theorem}[section]
\newtheorem{lemma}[thm]{Lemma}
\newtheorem{prop}[thm]{Proposition}

\theoremstyle{definition}
\newtheorem{defn}[thm]{Definition}
\newtheorem{remarka}[thm]{Remark}

\newtheorem{cor}{Corollary}[thm]

\theoremstyle{remark}
\newtheorem{remark}[thm]{Remark}

\newcommand{\R}{\mathbb{R}}

\newcommand{\C}{\mathbb{C}}
\newcommand{\Z}{\mathbb{Z}}
\newcommand{\D}{\mathbb{D}}

\newcommand{\bigO}{\mathcal{O}}
\newcommand{\s}{\mathbb{S}}

\newcommand{\ds}{\displaystyle}

\def\proof{\par\medskip\noindent {\bf Proof.\ \ }}

\def\qed{\hfill $\square$\\ }

\def \pvec #1#2{\begin{pmatrix}#1\\ \noalign{\vskip  -0 pt} #2\end{pmatrix}}
\def \hvec #1#2{\left(#1, #2\right)\!}

\def \diag #1#2#3#4#5#6#7#8{\begin{CD} #1 @>#5>> #2\\ @V#6VV  @VV#7V\\ #3 @>#8>> #4 \end{CD}}
\def\He{{H\'enon }}

\title[Complex H\'enon maps and discrete groups]{Complex H\'enon maps and discrete groups}
\author[Raluca Tanase]{Raluca Tanase$^1$}
\address{$^1$Institute for Mathematical Sciences, Stony Brook University, Stony Brook, NY 11794.}
\email{rtanase@math.sunysb.edu}

\date{\today}

\begin{document}

\maketitle

\vspace{.5cm}
\begin{abstract}
\noindent {\sc abstract.}
Consider the standard family of complex H\'enon maps
$H(x,y) = (p(x) - ay, x)$, where $p$ is a quadratic polynomial
and $a$ is a complex parameter. Let $U^{+}$ be the set of points that escape to
infinity under forward iterations. 
The analytic structure of the escaping set $U^{+}$ is well understood from previous
work of J. Hubbard and R. Oberste-Vorth as a quotient of 
$(\mathbb{C}-\overline{\mathbb{D}}) \times \mathbb{C}$ by a discrete
group of automorphisms $\Gamma$ isomorphic to
$\mathbb{Z}[1/2]/\mathbb{Z}$. On the other hand, the boundary $J^{+}$
of $U^{+}$ is a complicated fractal object on which the H\'enon map
behaves chaotically. We show how to extend the group action to
$\mathbb{S}^1\times\mathbb{C}$, in order to represent the set $J^{+}$
as a quotient of $\mathbb{S}^1\times \mathbb{C}/\Gamma$ by an
equivalence relation. We  analyze this extension for H\'enon
maps that are small perturbations of hyperbolic polynomials with
connected Julia sets or polynomials with a parabolic fixed point.
\end{abstract}
\vspace{.5cm}

\tableofcontents

\newpage 
\section{Introduction}\label{sec:Introduction}
\He maps have played an important role in the development of modern dynamics, both in the real and in the complex setting. Real \He maps were first introduced by Michel H\'enon as a simplified
model of the Poincar\'e section of the Lorenz model.  The dynamics of
\He maps is intriguing and challenging and they are some of the most
studied examples of dynamical systems that exhibit chaotic behavior. As a complex system, the \He map is also 
of major interest, due to the fact that all polynomial automorphisms of $\C^2$ can be reduced to
compositions of \He maps with simpler functions, as shown by S. Friedland and J. Milnor in \cite{FM}.

We  consider the standard family of complex H\'enon maps $H_{p,a}(x,y) = (p(x) - ay, x)$, where $p$ is a quadratic polynomial and $a$ is a complex parameter. Let $U^+$ and $U^-$ be the set of points that
escape to infinity under forward and respectively backward iterations of the H\'enon map. The topological boundaries $J^+$ of $U^+$ and $J^-$ of $U^-$ are complicated fractal sets on which the \He map
behaves chaotically. The sets $J^+$, $J^-$ and $J=J^+ \cap J^-$ are called the Julia sets of the \He map, and $J$ is considered to be the analogue of the Julia set from one-dimensional dynamics.

This article is devoted to discrete group actions
and connections with the topology of the set $J^+$. The analytic
structure of the escaping set $U^+$ is well understood from previous
work of J. Hubbard and R. Oberste-Vorth in \cite{HOV1} as a quotient
of $(\mathbb{C}-\overline{\mathbb{D}}) \times \mathbb{C}$ by a
discrete group of automorphisms $\Gamma$ isomorphic to
$\mathbb{Z}[1/2]/\mathbb{Z}$. As usual, $\mathbb{D}$ denotes the open unit disk in the complex plane. 
We explain this result in Section \ref{Covering}.

In Section \ref{sec: ExtendGroupJ+} we show how to extend the group
action to the boundary $\mathbb{S}^1\times\mathbb{C}$ in certain cases, in order to represent the
fractal set $J^+$ as a quotient of $\mathbb{S}^1\times
\mathbb{C}/\Gamma$ by an explicit equivalence relation. The group
extension has important topological consequences that we describe in Section
\ref{sec: Results}, where we analyze the extension for H\'enon maps
that are  perturbations of hyperbolic polynomials with connected Julia set. In Theorem \ref{thm:discont} we show that the group acts properly discontinuous and without fixed points on $\mathbb{S}^1\times\mathbb{C}$ and thus taking the quotient of $\mathbb{S}^1\times\mathbb{C}$ by the group action gives a topological manifold $\mathcal{M}$. The dynamics of the \He map on the set $J^+$ is semi-conjugate to the dynamics of a model map on $\mathcal{M}$. The semi-conjugacy function can be viewed as a two-dimensional analogue of the Carath\'eodory loop from polynomial dynamics. 
In the simplest case, when $p$ has an attractive fixed point ($p$ is taken from the interior of the main cardioid of the Mandelbrot set), an actual conjugacy is achieved, so $J^+$ itself is a topological manifold. In the other cases studied, we show that the set $J^+$ is a quotient of the manifold $\mathcal{M}$ by an equivalence relation which is described explicitly in Theorem \ref{thm: JpC}. 

The proof uses some results of M. Lyubich and J. Robertson 
\cite{LR} on the characterization of the critical locus for complex \He maps. The proof also
requires a careful analysis of the invariants of the \He map. In Section \ref{sec: ExtendGroupJ+} we
 introduce an important function for the study of the \He family, which we denote
$\alpha:\s^1\rightarrow \C^*$ and which encodes the dynamics of the \He map. The image of $\alpha$ is a fractal set, for which we have
designed and implemented a plotting algorithm in Section
\ref{Algorithm}. In Section \ref{Degeneracy}, we studied the degeneracy of the
cocycle $\alpha$ as the Jacobian tends to $0$, and came up with
an interesting relation connecting $\alpha$ with the group
action on $\s^1\times\C$. Section \ref{sec:grgamma} provides some
sharp estimates of the growth of the group elements. These are
useful for proving in Section \ref{sec:
Results} that the group 
acts properly discontinuously on $\s^1\times \C$.

\noindent \textit{Acknowledgements.} I would like to thank John Hubbard for
his entire support and guidance with this project.
I would also like to thank Xavier Buff and Remus Radu for 
many insightful conversations on this topic.

\section{Preliminaries}\label{sec:Preliminaries}

Consider the complex H\'enon map 
$H_{p,a}\hvec{x}{y} = \hvec{p(x)-ay}{x}\,$, 
where $p(x)=x^{2}+c$ is a monic polynomial of degree two. If $a\neq 0$, $H_{p,a}$ is a biholomorphism with constant Jacobian equal to
$a$, and the inverse map is 
$
H_{p,a}^{-1} \hvec{x}{y} = \hvec{y}{(p(y)-x)/a}.
$

The {\it filled-in Julia set}  of the polynomial $p$ is defined as
\[K_p = \{z\in \C\ :\ |p^{\circ n}(z)|\ \mbox{remains bounded as}\ n\rightarrow \infty \}.\]

The set  $J_p=\partial K_p$ is the {\it Julia set} of $p$. 
In analogy with one-dimensional dynamics, one defines the following dynamically invariant sets for the H\'enon map:
\[
K^{\pm}=\left\{ \hvec{x}{y}\in \C^2 : \left\| H_{p,a}^{\circ
n}\hvec{x}{y}\, \right\| \ \mbox{remains bounded as}\
n\rightarrow\pm\infty\right\}
\]
\[ U^{\pm} = \C^2-K^{\pm} \ \ \ \mbox{(the escaping sets)}
\]
\[
J^{\pm}=\partial K^{\pm}=\partial U^{\pm}
\]
\[
K=K^{-}\cap K^{+} \ \ \ \mbox{and}\ \ \ J=J^{-}\cap J^{+}.
\]

\noindent The sets $J^+$ and $J^-$ are closed, unbounded, connected
fractal objects in $\C^2$ \cite{BS1}. In the cases that we will be working
with, the Jacobian $a$ has absolute value less than $1$, so $K^{-}$ has no interior and 
$J^{-}=K^{-}$. When $a$ is small and $p(x)=x^2+c$ is a hyperbolic
polynomial, the interior of $K^{+}$ consists of the basins of
attraction of an attractive periodic orbit. The common boundary of
the basins is $J^+$ \cite{FS}, \cite{BS1}.

\vglue 0.25cm
\begin{figure}[htb]
\begin{center}
\includegraphics[scale =1.01]{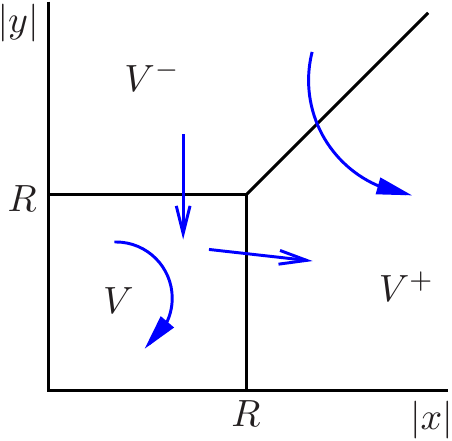}
\end{center}
\caption{Filtration of $\C^2$.}
\label{fig:filtration}
\end{figure}

\noindent According to \cite{HOV1}, for $R>2$ sufficiently large, the dynamical space $\C^2$ can be divided
into three regions:
$ V=\{(x,y)\in \C^2: |x|\leq R, |y|\leq R\},$
\[
V^+=\{(x,y):|x|\geq \max(|y|,R) \}\ \
\mbox{and}\ \ V^-=\{(x,y): |y|\geq \max(|x|, R)\}.
\]
\noindent The sets $J$ and $K$ are contained in the polydisk $V$. The escaping sets can be
described as union of backward iterates of $V^+$ and respectively forward iterates of $V^-$ under the \He map:
\[
U^+=\bigcup_{k\geq 0}
H^{-\circ k}(V^+) \ \ \mbox{and}\  \  U^-=\bigcup_{k\geq 0} H^{\circ k}(V^-).
\]
\noindent The domains $V^{+}$ and $V^{-}$ are easier to understand
because one can define an analogue of the B\"ottcher
coordinates. More precisely we have the following lemma:

\begin{lemma}[Hubbard, Oberste-Vorth \cite{HOV1}]\label{lemma:varphi}
There exists a unique holomorphic map
$\varphi^+:V^+\rightarrow \C-\overline{\D}$ such that
$
\varphi^+\circ H = (\varphi^+)^2
$
and $\varphi^+(x,y)\sim x$ as $(x,y)\rightarrow\infty$ in $V^+$. There exists a unique holomorphic map
$\varphi^-:V^-\rightarrow \C-\overline{\D}$ such that
$
a\varphi^-\circ H^{-1} = (\varphi^-)^2
$
and $\varphi^-(x,y)\sim y$ as $(x,y)\rightarrow\infty$ in $V^-$.
\end{lemma}

The set $U^+$ is foliated by copies of $\C$, which have a natural affine
structure. The holomorphic function $\varphi^+$ defines a holomorphic foliation on $V^+$.
The leaves of the foliation are just the level sets of $\varphi^{+}$. One can then extend this foliation from $V^{+}$ to $U^{+}$ by the dynamics.
The function ${(\varphi^+)}^{2^k}$ is well defined on $H^{-\circ k}(V^+)$ as
${(\varphi^+)}^{2^k}=\varphi^+\circ H^{\circ k}$ and it defines a holomorphic foliation on $H^{-\circ k}(V^+)$.

\vglue 0.5cm
\begin{figure}[htb]\label{fig:F2}
\begin{center}
\hspace{-1cm}
\includegraphics[scale =0.67]{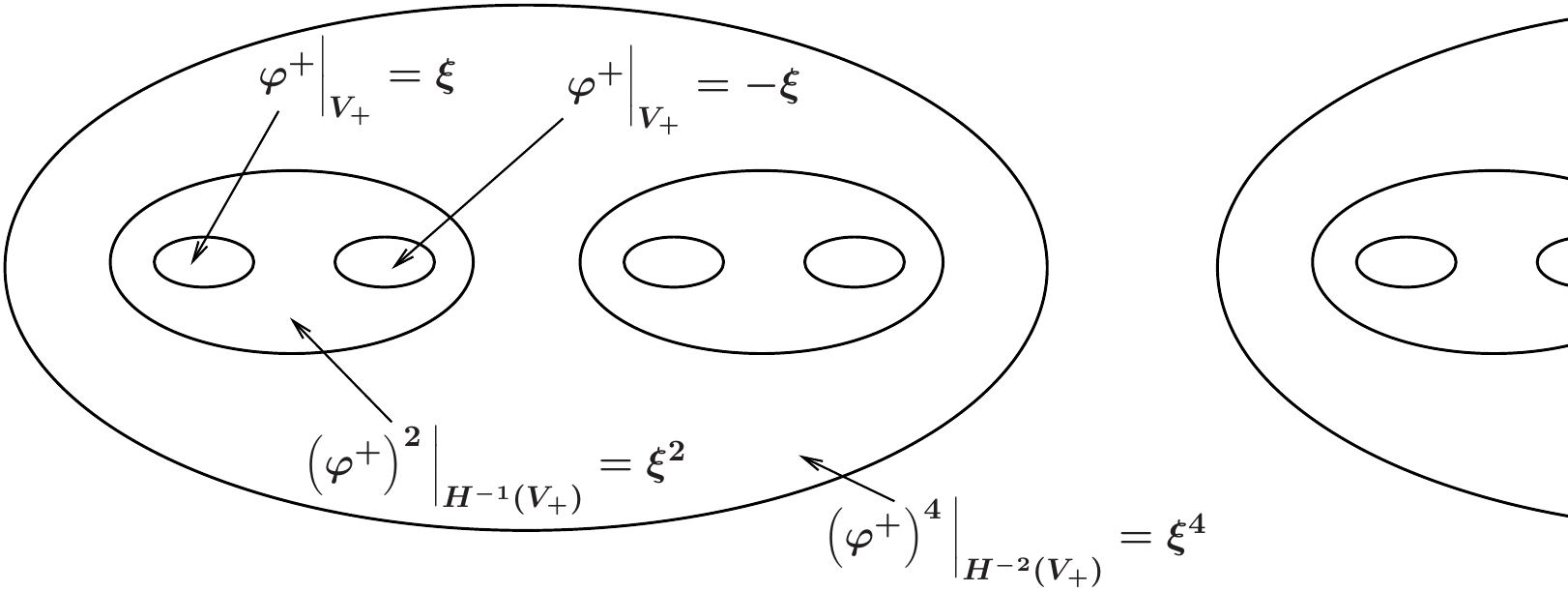}
\end{center}
\caption{A fiber $\mathcal{F}_{\xi}$ of the foliation of $U^+$, for $\xi\in \C-\overline{\D}$.}
\end{figure}

\noindent One can also define a similar holomorphic foliation of the set $U^{-}$ using the map $\varphi^{-}$.
The foliations of the escaping sets $U^+$ and $U^-$ are not everywhere transverse to each other. The critical locus $\mathcal{C}$ of the \He map is the set of
tangencies between the foliation of $U^+$ and the foliation of
$U^-$. The set $\mathcal{C}$ is a closed analytic subvariety of
$U^+\cap U^-$ and is invariant under the \He map.
\begin{thm}[Bedford, Smillie \cite{BS5}]
The critical locus $\mathcal{C}$ is nonempty. The boundary
$\partial{\mathcal{C}}$ of $\mathcal{C}$ intersects both $J^+$ and
$J^-$ and we have $\overline{\mathcal{C}}\cap J^+\cap U^-\neq
\emptyset$ and $\overline{\mathcal{C}}\cap J^-\cap U^+\neq
\emptyset$.
\end{thm}
\begin{thm}[Lyubich, Robertson  \cite{LR}]\label{thm: LyubichRobertson}
\justifying Let $H$ be a hyperbolic \He map with connected $J$,
which is a small perturbation of a hyperbolic polynomial
$p(x)=x^2+c$, with connected Julia set $J_p$. We have the following
description of the critical locus:
\begin{itemize}
  \item[(a)]  \justifying There exists a unique primary component $\mathcal{C}_0$ of the {\it
critical locus} asymptotic to the $x$-axis.
  \item[(b)] \justifying There exists a biholomorphic extension of $\varphi^+$ from
$\mathcal{C}_0$ to $\C-\overline{\D}$.
  \item[(c)] \justifying Moreover, there exists a biholomorphism $\tau^+$  from
$\mathcal{C}_0$ to $\C-K_p$, which can be extended homeomorphically
from $\overline{\mathcal{C}}_0$ to $\C-\mathring{K_p}$.
   \item[(d)] \justifying $\mathcal{C}_0$ is everywhere transverse to the
   foliation of $U^+$ and $U^-$.
   \item[(e)] \justifying All other components of $\mathcal{C}$ are
forward or backward
   iterates of $\mathcal{C}_0$ under $H$.
\end{itemize}
\end{thm}

\begin{remarka}
Since the \He map is hyperbolic with connected Julia set, the
boundary of $\mathcal{C}_0$ belongs to $J^+$. The forward iterates
of $\mathcal{C}_0$ accumulate on $J^-$.
\end{remarka}
\begin{remarka} A model for the critical locus is also described in \cite{F} for perturbations of quadratic hyperbolic polynomials with disconnected Julia sets. The critical locus is connected in this case.
\end{remarka}
\noindent{\bf The degenerate case $a=0$.} The picture when $a$ is $0$ helps visualize the foliation of $U^+$ and $J^+$ and the primary component of the critical locus. The \He map $ H_{p,0} (x,y)=\hvec{p(x)}{x}\,$
is no longer a biholomorphism and maps all $\C^2$ to the curve
$\{x=p(y)\}$. However, the foliations of $U^+$ and $J^+$ persist and are easier to describe:
\begin{itemize}
\item[(a)] $\varphi^+$ is just the B\"ottcher
isomorphism of $p$.
\item[(b)] $J^+ = J_p\times \C$, where $J_p$ is the Julia set of $p$.
\item[(c)] $U^+ = (\C - K_p)\times \C$, where $K_p$ is the
filled-in Julia set of $p$.
\end{itemize}
\noindent The primary component
of the critical locus can also be easily understood from \cite{LR},
$\overline{\mathcal{C}}_0=(\C - \mathring{K_p})\times \C$,
where $\mathring{K_p}$ is the interior of the filled-in Julia
set of $p$.

\section{The covering space of the escaping set $U^{+}$}\label{Covering}
In this section we describe the analytic
structure of the escaping set $U^+$. 
\begin{lemma}
There exists a closed holomorphic $1$-form on $U^+$, with $H^*w=2w$.
\end{lemma}
\proof The map $\varphi^+$ is well defined on $U^+$ up to local choices of roots of
unity, so $\log(\varphi^+)$ is well defined up to local addition of
constants. Hence the form $w=d\log\varphi^+$ is well defined and holomorphic on $U^+$.
It is easy to see from Lemma \ref{lemma:varphi}
that 
\[
H^*w=H^*\frac{d\varphi^+}{\varphi^+}=\frac{d(\varphi^+\circ
H)}{\varphi^+\circ H}=\frac{d((\varphi^+)^2)}{\varphi^+}=2w.
\]
\vglue -0.7cm
\qed

\begin{defn}
For a closed curve $C$ in $U^+$, define the index $\eta(C )$ as
\begin{equation}\label{eq:etaC}
\eta(C ):=\frac{1}{2\pi i}\int_{C}w.
\end{equation}
\end{defn}

\noindent Since $w$ is a closed 1-form, the number $\eta(C )$ depends only on the homotopy type of $C$.
The following properties from \cite{BS8} and \cite{MNTU} of $\eta$ are helpful for understanding the topology of the escaping set $U^{+}$:
\begin{itemize}
  \item[(a)]
  $\ds \eta(H^k(C))=\frac{1}{2\pi i}\int_{H^k(C)}w=\frac{1}{2\pi i}\int_{C}(H^k)^*w=2^k\eta(C)$.
  \item[(b)] Take $C\subset V^+$. One can homotopically enlarge $C$ so
  that it belongs to the region where $\varphi^+(x,y)\sim x$. Then
  $ \eta(C)=\frac{1}{2\pi i}\int_{C}\frac{dx}{x}\in \Z$, since
  it represents the winding number of $C$ around the $x$-axis.
  \item[(c)] In $U^+$, there exists $k$ such that $H^k(C)\subset V^+$, so $\eta(H^k(C))=m\in
  \Z$.
  Therefore $\eta(C)=\frac{m}{2^k}$, where $m$ and $k$ are integer numbers.
\end{itemize}

\begin{lemma}[\cite{HOV1}]\label{lemma:Z1/2}
The fundamental group of $U^{+}$ is isomorphic to $\Z[\frac{1}{2}]$ where $\Z[\frac{1}{2}]=\left\{\frac{m}{2^{k}}\ |\ m,k\in \Z\right\}$.
\end{lemma}
\proof The proof is an immediate consequence of the properties (a), (b) and (c) listed above.
\qed

We would like of course to be able to extend $\varphi^{+}$, and not only $d\log \varphi^{+}$
to the whole set $U^{+}$. However there are topological obstructions which become
clear once we look at the behavior of the H\'enon map near infinity in $V^{+}$.
\begin{lemma}[\cite{HOV1}] For large $r$, the set $V^+(r)=\{(x,y)\in V^+ ,
|\varphi^+|=r\}$ is homeomorphic to a solid torus, and
$\varphi^+:V^+(r)\rightarrow \{z, |z|=r\}
$
is a fibration with fibers homeomorphic to closed disks.
On $V^+$ the H\'enon map is solenoidal, and the following diagram commutes:
\begin{figure}[htb]
\begin{center}
\includegraphics[scale =.775]{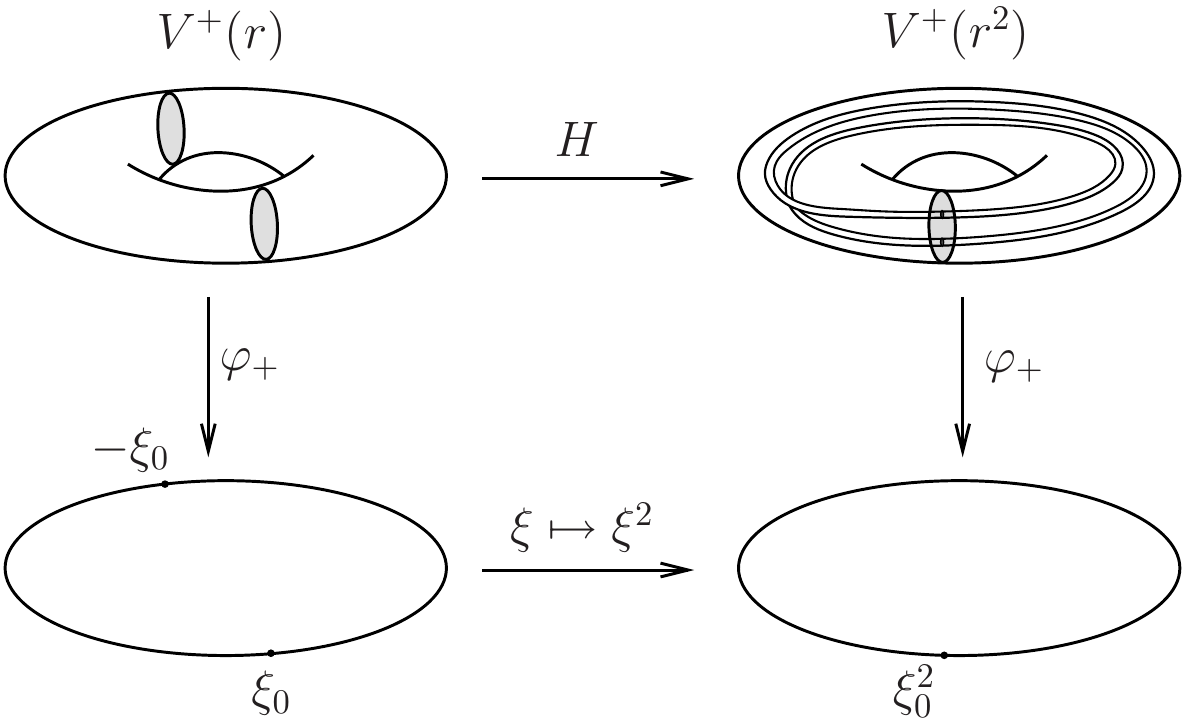}
\end{center}
\end{figure}
\end{lemma}

The map $\varphi^+$ does not extend holomorphically to $U^+$,
but it does extend along curves contained in $U^+$ which start in $V^+$
so it is well defined on a covering manifold $\widetilde{U}^+$ of
$U^+$. 
The covering manifold $\widetilde{U}^+$ is called the Riemann surface of $\varphi^+$ and its construction is fairly standard, nonetheless, for completion, we will outline  
an explicit construction of $\widetilde{U}^+$   from \cite{MNTU}. 

One starts by fixing a base point $a\in
V^+$, and defines the set $\widetilde{U}^+$ as follows:
\begin{eqnarray*}
\widetilde{U}^+&=&\left\{(z,C),\ z\in U^+,\ C\ \mbox{is a path in}\ U^+ \
\mbox{between}\ a\ \mbox{and}\ z\right\}\big{/} \sim\\
&&\ \mbox{where}\ (z,C)\sim(z',C')\Leftrightarrow z=z'\ \mbox{and}\
\eta(CC'^{-1})\in \Z.
\end{eqnarray*}
The definition does not depend on the choice of a particular
representative $(z,C)$ for the equivalence class. If $(z',C')$ is
another representative, then $z=z'$ and $\eta(CC'^{-1})=m$ $\in\Z$. It
follows that $\int_{C}w-\int_{C'}w=2\pi i m$, so
$\eta([z,C])=\eta([z',C'])$. There is also an analogue of $V^+$ in
$\widetilde{V}^+$, represented by the set
\[
    \widetilde{V}^+=\left\{[z,C]\in \widetilde{U}^+,\ z\in V^+,\
C\subset V^+\right\}.
\]
One can then define a lift of $\varphi^{+}$ to the covering manifold
$\widetilde{U}^+$. Let
$\widetilde{\varphi}^+:\widetilde{U}^+\rightarrow \C-\D$ be given by
the relation
\[
\widetilde{\varphi}^+([z,C])=\varphi^+(a)\exp\int_{C}w.
\]

\noindent It is easy to show that 
\begin{equation}\label{eq:02}
\widetilde{\varphi}^+\big{|}_{\widetilde{V}^+}=\varphi^+\big{|}_{V^+}\  \mbox{ and } \ 
\widetilde{U}^+ = \bigcup\limits_{n \geq 0}\widetilde{H}^{-\circ
n}(\widetilde{V}^+).
\end{equation} 
To check the first equality from Relation \ref{eq:02}, take an equivalence class $[z,C]\in \widetilde{V}^+$ such
that $z\in V^+$. One can verify by direct computation that
\[
\widetilde{\varphi}^+([z,C])=\varphi^+(a)\mbox{e}^{\log(\varphi^+(z))-\log(\varphi^+(a))}=\varphi^+(z).
\]
To show the second part of Relation \ref{eq:02}, take an equivalence class $[z,C]\in \widetilde{H}^{-\circ
n}(\widetilde{V}^+)$ such that $z\in H^{-\circ n}(V^+)$. Then we have
\[
{\widetilde{\varphi}^+}([z,C])=\varphi^+(a)\exp\int_{C}w=\varphi^+(a)\exp
\left (\frac{1}{2^n}\int_{H^{\circ n}C}w\right )
\]
and it follows that
\begin{eqnarray}\label{eq:03}
\ds\widetilde{\varphi}^{+}([z,C])^{2^n}&=&\varphi^+(a)^{2^n}\mbox{e}^{\log(\varphi^+(H^{\circ
n}(z)))-\log(\varphi^+(H^{\circ
n}(a)))} \nonumber \\
&=&\varphi^+(H^{\circ n}(z))=\varphi^+(z)^{2^n}.
\end{eqnarray}

\begin{thm}[Hubbard, Oberste-Vorth \cite{HOV1}] \label{thm:HOV}The covering manifold  $\ \widetilde{U}^+$ is a trivial
analytic fiber bundle over $\C-\overline{\D}$, with fibers
isomorphic to $\C$.
\end{thm}

A nice proof of this theorem is given in \cite{HOV1} and we will not reproduce it here in detail. The key point of the proof is to show that  the map $\widetilde{\varphi}^+:\widetilde{U}^+\rightarrow
(\C- \overline{\D})$ is an analytic submersion with fibers
isomorphic to $\C$, then show (by a nontrivial argument) that $\widetilde{U}^+$ is a locally trivial fiber bundle, locally homeomorphic to $(\C- \overline{\D})\times \C$. Then the result of Theorem \ref{thm:HOV} follows by complex analysis, as $\C-
\overline{\D}$ is Stein, so topological and analytic classification
of line bundles over $\C- \overline{\D}$ coincide.

\tikzset{node distance=3cm, auto} 
\begin{center}
\begin{tikzpicture}
  \node (C) {$\widetilde{U}^+$};
  \node (B) [below of=C] {$U^{+}\supset V^{+}$};
  \node (A) [right of=B] {$\C-\overline{\D}$};
  \draw[decoration={markings,mark=at position 1 with {\arrow[scale=2]{>}}},
    postaction={decorate}] (C) to node {$\widetilde{\varphi}^{+}$} (A);
  \draw[decoration={markings,mark=at position 1 with {\arrow[scale=2]{>}}}, 
    postaction={decorate}] (C) to node [swap] {$\pi$} (B);
  \draw[decoration={markings,mark=at position 1 with {\arrow[scale=2]{>}}},
    postaction={decorate}] (B) to node [swap] {$\varphi^{+}$} (A);
\end{tikzpicture}
\end{center}


\begin{thm}[Hubbard, Oberste-Vorth  \cite{HOV1}]\label{thm: HOV-escapingSet}
The analytic structure of $U^+$ is well-understood:
\begin{equation*}
U^{+} = (\C- \overline{\D})\times \C \big{/}\Gamma_{p,a}, 
\end{equation*}
where $\Gamma_{p,a}\subset\mbox{Aut}((\C- \overline{\D})\times \C)$
is a discrete group isomorphic to $\Z[\frac{1}{2}]\big{/}\Z$.
\end{thm}
\proof By Theorem \ref{thm:HOV},  $\widetilde{U}^+$ is a covering manifold of $U^+$, hence one can describe $U^+$ as a quotient of
$\widetilde{U}^+$ by a group $\Gamma_{p,a}$ of deck transformations. The group $\Gamma_{p,a}$ is isomorphic
to $\pi_1(U^+)\big{/}\pi_1(\widetilde{U}^+)$, hence  isomorphic to
$\Z\left[\frac{1}{2}\right]\big{/}\Z$ by Lemma \ref{lemma:Z1/2}. The following diagram depicts the situation:
\begin{equation*}
\diag{\widetilde{U}^+}{(\C-\overline{\D})\times \C}{U^+}{(\C-
\overline{\D})\times \C
\big{/}\Gamma_{p,a}}{\simeq}{\pi}{\pi}{\simeq}
\end{equation*}
\vglue -0.45cm
\qed

There is a unique lift $\widetilde{H}$ of $H$ to the covering
manifold $\widetilde{U}^{+}$ such that the following diagram
commutes and conditions 1-4 hold:

\hspace{1cm}
\begin{minipage}[t]{0.17\textwidth}
\vglue -0.2cm
\[
\diag{\widetilde{U}^+}{\widetilde{U}^+}{U^+}{U^+}{
\widetilde{H}}{\pi}{\pi}{H}
\]
\end{minipage}
\hspace{1.5cm}
\begin{minipage}[t]{0.75\textwidth}
\vglue 0.02cm
\begin{itemize}
\item[1.] $\pi\circ  \widetilde{H} = H\circ \pi$
\item[2.] $\widetilde{\varphi}^{+}\circ  \widetilde{H}= \left(\widetilde{\varphi}^{+}\right)^{2}$ on
$\widetilde{U}^{+}$
\item[3.] $\widetilde{H}\circ \gamma = (2\gamma)\circ \widetilde{H}$, for all $\gamma\in
\Gamma_{p,a}$
\item[4.] $\widetilde{H}(\widetilde{V}^+)\subset \widetilde{V}^+ $
\end{itemize}
\end{minipage}

\vspace{0.5cm} 
The map $\widetilde{H}$ is a covering map from
$\widetilde{U}^+$ to $\widetilde{U}^+$ with sheet number $2$.

\begin{remarka} The foliation of the covering manifold
$\widetilde{U}$ by level sets of the function
$\widetilde{\varphi}^+$ descends to a foliation of the escaping set
$U^+$. This is the same foliation as the  one induced by the function $\varphi^+$ on $U^+$
in Section \ref{sec:Preliminaries}, as it can be easily seen from relations \ref{eq:02} and \ref{eq:03} and the properties of
the lift $\widetilde{H}$.
\end{remarka}

\section{The Stable Multiplier Condition}\label{sec: ExtendGroupJ+}
We will show how to extend the group action
$\Gamma_{p,a}$ to $\s^1\times\C$ in certain cases, in order to
represent the fractal boundary $J^+$ of $U^+$ as a quotient of
$\s^1\times \C/\Gamma_{p,a}$ by an equivalence relation. We will
discuss this in Theorem \ref{thm: main1} and Corollary \ref{cor:Mc}.

Let us first explain the meaning of an extension of the group
elements to $\s^1\times \C$. After a particular trivialization of
the covering manifold $\widetilde{U}^+$ has been chosen, one can
define a lift of the \He map to $(\C-\overline{\D})\times \C$ so
that the following diagram commutes
\[
\diag{(\C-\overline{\D})\times \C}{(\C-\overline{\D})\times
\C}{U^+}{U^+}{\widetilde{H}}{\pi}{\pi}{H}
\]

\noindent In follows from conditions 1-3 that the lift
$\widetilde{H}$ of the \He map is an analytic function of the form
\begin{equation}\label{eq:tildeH}
\widetilde{H}\hvec{\xi}{z}=\hvec{\xi^2}{\alpha(\xi)z+\beta(\xi)},
\end{equation}
where $\alpha:\C-\overline{\D}\rightarrow \C^*$ and 
$\beta:\C-\overline{\D}\rightarrow \C$ are analytic functions.

\medskip
\noindent {\bf Extension.} We would like to
 extend the map $\widetilde{H}$ to $\s^1\times \C$, so that the dynamics of
$\widetilde{H}$ on $\s^1\times \C$ is "compatible" with the dynamics
of $H$ on the Julia set $J^+$. However, the set $J^+$ contains stable
manifolds of periodic points in $J$ (as we will see in Theorem \ref{thm:
BS-lamination JU}), so we must require that the
following condition is satisfied:

\medskip
\noindent {\bf Stable Multiplier Condition.} The functions $\alpha$ and $\beta$ extend
continuously to $\s^{1}$. The stable multipliers of $\widetilde{H}$
on $\s^1\times\C$ agree with the stable multipliers of $H$ on $J^+$,
in the sense that for every periodic point $\xi=\xi^{2^k}$ of the
doubling map on $\s^1$ there exists a $k$-periodic point $x$ of
the \He map $H$ such that
\begin{equation}\label{eq:restriction}
    \alpha(\xi)\alpha(\xi^2)\ldots\alpha(\xi^{2^{k-1}})=\lambda(DH^{\circ
    k}_{x}),
\end{equation}
where $\lambda$ is the small eigenvalue of $DH^{\circ k}$ at $x$.

\medskip
We will call a function $\alpha$ which satisfies the Stable Multiplier Condition a cocycle.

\begin{remarka} In \cite{HOV1}, in addition to the description of the 
covering manifold $\widetilde{U}^+$,
it is also shown that there exists a unique isomorphism
$\widetilde{U}^+\rightarrow (\C-\overline{\D})\times \C$ such that,
with this trivialization, the map $\widetilde{H}$ is written as
\[
    \widetilde{H}\hvec{\xi}{z}=\hvec{\xi^{2}} {\frac{a}{2}z +\xi^{3}-\frac{c}{2}\xi}.
\]
Notice that even if there is no problem in continuously
extending this map to $\s^1\times \C$, the dynamics of the extension
$\widetilde{H}$ on $\s^1\times\C$ is quite different from the
dynamics of the \He map $H$ on $J^+$. The stable multipliers of
$\widetilde{H}$ are ``too simple'', as they are all equal to $a/2$,
whereas the multipliers of the \He map can (and will) be
complicated.
\end{remarka}

In Section \ref{sec:Preliminaries}, we described the foliation of the escaping set $U^+$. When the \He map $H$ is hyperbolic, the boundary $J^+$ of $U^+$ is also foliated by copies of $\C$, given by stable manifolds of points from the Julia set $J$, as illustrated by the following theorem: 

\begin{thm}[Bedford, Smillie \cite{BS7}]\label{thm: BS-lamination JU}\hspace{2cm}
\begin{itemize}
  \item[(a)] Let $p\in J$ be a saddle periodic point of the \He map, then $J^+$ is the closure of
$W^s(p)$.
  \item[(b)] For hyperbolic \He maps with Jacobian $|a|<1$, the set
$J^{+}=W^s(J)$, so $J^+$ has its own dynamically defined Riemann
surface lamination, whose leaves consist of the stable manifolds
$W^s(p)$ of points $p\in J$.
  \item[(c)]  If in addition, the Julia set $J$ is connected, then the foliation
of $U^+$ and the lamination of $J^+$ fit together continuously to
form a locally trivial lamination of $U^+\cup J^+$.
\end{itemize}
\end{thm}
\noindent When $H$ is hyperbolic, for each point $p\in J$ there
exists a biholomorphic function $\varphi:\C\rightarrow W^s(p)$ which
 defines an affine structure on the stable manifold
$W^s(p)$. In addition, the iterates of the \He map $H$ preserve the
affine structure, in the sense that the pull-back or push-forward of
the affine structure from one leaf to another agrees with the
original affine structure on the new leaf \cite{BS5}.

\section{A trivialization of the lamination of $U^+\cup J^+$} \label{Affine}

We will use the primary component of the critical locus from Theorem \ref{thm: LyubichRobertson} to give an
identification of each of the fibers $\mathcal{F}_{\xi}=(\widetilde{\varphi}^+)^{-1}(\xi)$ with $\C$.
These identifications will provide a specific trivialization of the
bundle $\widetilde{U}^+ \simeq (\C-\D)\times \C$. 

It follows from Theorem \ref{thm: LyubichRobertson}
that the primary component $\mathcal{C}_0$ of the critical locus is biholomorphic to the
exterior of the filled-in Julia set $K_p$ of the
polynomial $p$ via a map $\tau^{+}: \mathcal{C}_0\rightarrow \C- K_p$  that extends to a homeomorphism between the boundary of $\mathcal{C}_0$ and the Julia set $J_p$. Therefore, the closure of the primary component
$\overline{\mathcal{C}}_0$ can be naturally identified with $\C-\D$
via the map
\begin{equation*}
\C-\D\xrightarrow[\ \ \ \ \ \
]{\gamma}\C-\overset{\circ}{K_p}\xrightarrow{(\tau^{+})^{-1}}\overline{\mathcal{C}}_0,
\end{equation*}
where the composition
$\left(\tau^{+}|_{\overline{\mathcal{C}}_0}\right)^{-1}\circ \gamma$
is biholomorphic on $\C-\overline{\D}$ and continuous on $\C-\D$. The
function $\gamma$ in the diagram above is the B\"ottcher  coordinate
of the polynomial $p$.

We briefly recall the definition of the B\"ottcher  coordinate from
\cite{DH} and \cite{M}. Let $p$ be a quadratic polynomial with connected filled-in Julia set
$K_p$. There exists a unique analytic map $\varphi:
\C-K_p\rightarrow \C-\overline{\D}$ tangent to the identity at
infinity that conjugates $p$ to $z\rightarrow z^2$, that is
$\varphi\circ p =(\varphi)^2$. The function $\varphi$ is called the
\emph{B\"ottcher isomorphism}, and the inverse map
$\gamma=\varphi^{-1}:\C-\overline{\D}\rightarrow \C-K_p$ the
\emph{B\"ottcher coordinate}. If in addition the filled-in Julia set $K_p$ is locally connected,
the B\"ottcher  coordinate extends continuously to $\gamma:\s^1\rightarrow J_p$,
$t\rightarrow \lim\limits_{r\rightarrow 1^+} \gamma(re^{2 \pi i t})$. The extension is a continuous surjective map called the \emph{Carath\'{e}odory loop}.

\begin{lemma}[\textbf{Trivialization lemma}]\label{lemma: pi} There exists a continuous surjective function
$\pi:(\C-\D)\times \C \rightarrow U^+\cup J^+$, holomorphic from
$(\C-\overline{\D})\times \C \rightarrow U^+$ and analytic on the
leaves of the lamination of $J^+$, such that the following diagram
commutes
\[
\diag{(\C-\D)\times \C}{(\C-\D)\times \C}{U^+\cup J^+}{U^+\cup
J^+}{\widetilde{H}}{\pi}{\pi}{H}
\]
where $\widetilde{H}\hvec{\xi}{z}=\hvec{\xi^2}{\alpha(\xi)z+\beta(\xi)}$, and
the functions $\alpha:\C-\D\rightarrow \C^{*}$ and
$\beta:\C-\D\rightarrow \C$ are continuous on $\C-\D$ and analytic
on $\C-\overline{\D}$.
\end{lemma}
\proof As in \cite{BV}, we will construct holomorphic parametrizations of the leaves of the foliation of $U^+$ that converge locally uniformly to the parametrization of a limit leaf of the lamination of $J^+$. Let $\mathcal{F}_{\xi}$ be a leaf of the lamination of
$U^+\cup J^+$. The critical points $c_0(\xi)$ and $c_{-1}(\xi)$
belong to $\mathcal{F}_{\xi}$ and they are given by the relation
\begin{equation}\label{eq:c0c1}
c_0(\xi)= \left(\tau^{+}|_{\overline{\mathcal{C}}_0}\right )^{-1}\circ \gamma(\xi)\ \  \mbox{ and } \ \
c_{-1}(\xi)=H^{-1}(c_0(\xi^2)).
\end{equation}
\noindent Each leaf is biholomorphic to $\C$ and there exists a unique
analytic mapping\\ $\pi_{\xi}:\C\rightarrow\mathcal{F}_{\xi}$ which sends
\begin{equation}\label{eq:sttriv}
    0\rightarrow c_0(\xi)\ \  \mbox{ and } \ \ 1\rightarrow c_{-1}(\xi).
\end{equation}
We can therefore define the function $\pi: (\C-\D)\times
\C\rightarrow U^+\cup J^+$ by $\pi(\xi,z)=\pi_{\xi}(z)$.

Recall from Section \ref{sec:Preliminaries} that $\mathcal{F}_{\xi}=\mathcal{F}_{\omega\xi}$
for all $\omega$ dyadic roots of unity ($w$ is dyadic if $w^{2^k}=1$, for some
non-negative integer $k$). Of course, the primary component $\mathcal{\overline{C}}_0$  of the
critical locus intersects
$\mathcal{F}_{\xi}$ at all points of the form $c_0(\omega\xi)$,
where $\omega^{2^k}=1$ for some integer  $k\geq 0$. Therefore we
will end up parametrizing the same fiber $\mathcal{F}_{\xi}$ "a
dyadic number of times". We parametrize $\mathcal{F}_{\omega\xi}$
by first fixing $c_0(\omega\xi)$ at the origin. The
nearest intersection point of $\overline{\mathcal{C}}_{-1}$ with
$\mathcal{F}_{\xi}$ will then be $c_{-1}(\omega\xi)$, and we set
this to be $1$.

The function $c_0:\C-\D\rightarrow \C^2$ is holomorphic on
$\C-\overline{\D}$ and continuous on $\s^1$, as shown in Theorem
\ref{thm: LyubichRobertson}. Consequently $c_{-1}:\C-\D\rightarrow
\C^2$ has the same properties. The primary component of the critical
locus $\overline{\mathcal{C}}_{0}$ is transverse to the leaves on
the foliation of $U^+$, by Theorem \ref{thm: LyubichRobertson}. The
affine structure on $J^{+}\cup U^{+}$ is transversely continuous
\cite{BS5}. Hence the projection $\pi$ is continuous on
$(\C-\D)\times \C$.

The fact that the function $\pi$ is analytic on
$(\C-\overline{\D})\times \C$ follows from the construction of the
covering manifold $\widetilde{U}^+$. It is worth noting at this stage that
$\pi:(\C-\overline{\D})\times \C \rightarrow U^+$ is a covering map,
but $\pi: \s^1\times \C\rightarrow J^+$ is not in general a covering
map, unless the \He map is a perturbation of a quadratic polynomial
with an attractive fixed point.

The \He map becomes
$\widetilde{H}(\xi,z)=\left(\xi^2,\alpha(\xi)z+\beta(\xi)\right)$. For a fixed
$\xi$, we  compute $\alpha(\xi)$ and $\beta(\xi)$ by looking at
the affine structures on the fibers $\mathcal{F}_{\xi}$ and
$\mathcal{F}_{\xi^2}$.
\begin{eqnarray}\label{eq:system}
\alpha(\xi)[c_{-1}(\xi)]_{\mathcal{F}_{\xi}}+\beta(\xi) &=& [c_0(\xi^2)]_{\mathcal{F}_{\xi^2}}\\
\alpha(\xi)[c_{0}(\xi)]_{\mathcal{F}_{\xi}}+\beta(\xi) &=&
[H(c_0(\xi))]_{\mathcal{F}_{\xi^2}}\nonumber
\end{eqnarray}

\begin{figure}[htb]
\begin{center}
\includegraphics[scale =.75]{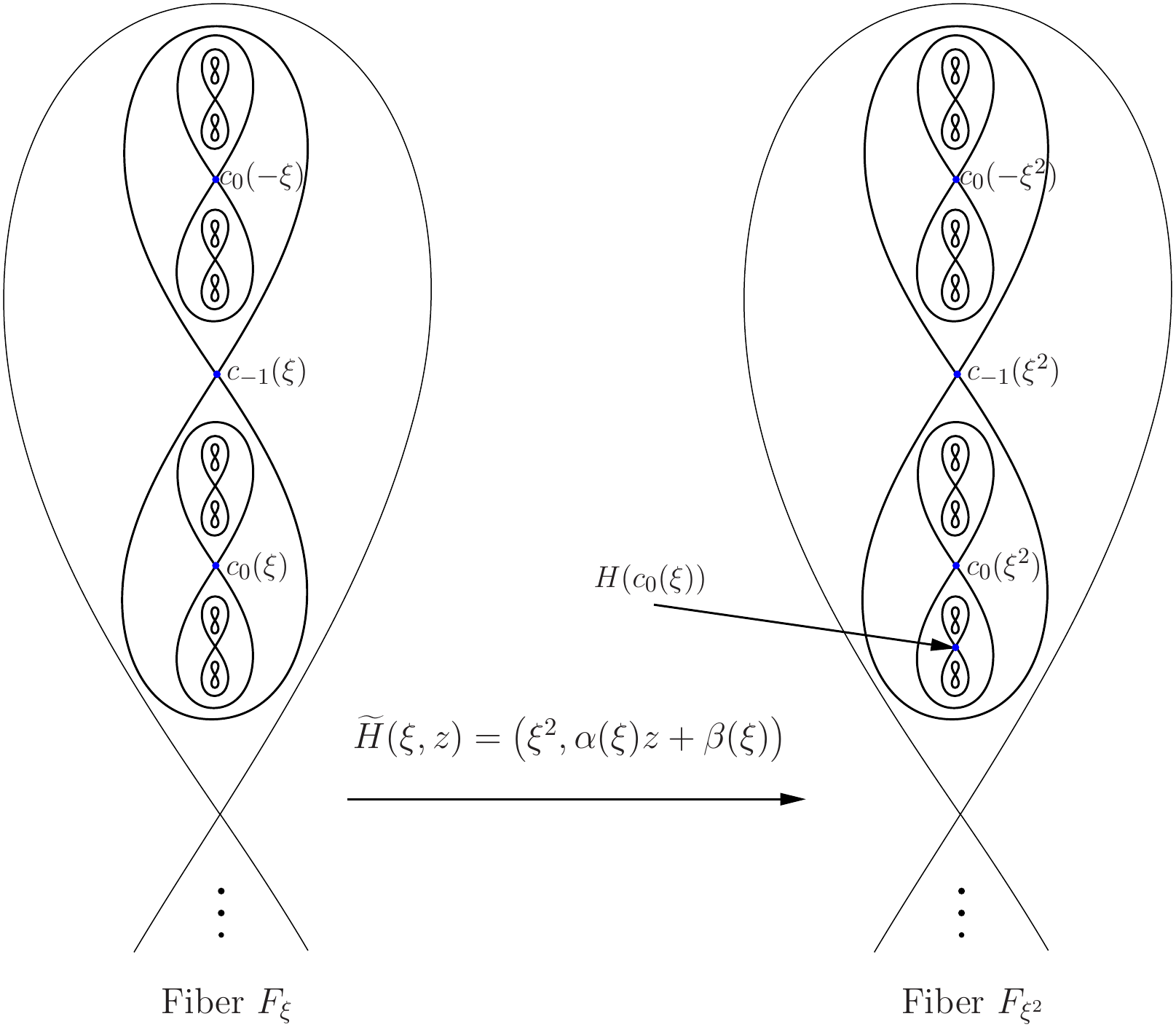}
\end{center}
\caption{Two fibers $\mathcal{F}_{\xi}$ and $\mathcal{F}_{\xi^{2}}$ of the lamination of $J^{+}\cup U^{+}$ and the action of the map $\widetilde{H}$ on the critical points $c_{0}(\xi)$ and $c_{-1}(\xi)$.} \label{pic:Fxi}
\end{figure}

The fiber $\mathcal{F}_{\xi^2}$ is biholomorphic to $\C$,
hence the ratio of the points $c_0(\xi^2)$, $H(c_0(\xi))$ and
$c_{-1}(\xi^2)$ does not depend on the choice of affine maps on
$\mathcal{F}_{\xi^2}$. Hence if we denote by
$[H(c_0(\xi))]_{\mathcal{F}_{\xi^2}}$ the coordinate of the point
$H(c_0(\xi))$ with respect to the particular affine map on
$\mathcal{F}_{\xi^2}$ which assigns $c_0(\xi^2)\rightarrow 0$ and
$c_{-1}(\xi^2)\rightarrow 1$, we obtain
\[
[H(c_0(\xi))]_{\mathcal{F}_{\xi^2}}=-\frac{H(c_0(\xi))-c_0(\xi^2)}{c_0(\xi^2)-c_{-1}(\xi^2)}
\]
and we can compute
\begin{equation}\label{eq:alpha-01}
\alpha(\xi)=-\beta(\xi)=\frac{H(c_0(\xi))-c_0(\xi^2)}{c_0(\xi^2)-H^{-1}(c_0(\xi^4))}.
\end{equation}
Clearly $\alpha$ does not vanish. Otherwise
$\overline{\mathcal{C}}_{0}$ and $\overline{\mathcal{C}}_{1}$ would
intersect and this is not possible from \cite{LR}. More precisely,
the primary component $\overline{\mathcal{C}}_{0}$ is inside a
trapping region around the $x-$axis that contains no other
components of the critical locus. \qed

\begin{remarka} We would like to write $c_1(\xi^2)$ in place of $H(c_0(\xi))$ in Equation \ref{eq:system},
but this would be incorrect, as we would not be able to distinguish
between $H(c_0(\xi))$ and $H(c_0(-\xi))$, which are two distinct
points of $\mathcal{F}_{\xi^2}$ (see also  Figure \ref{pic:Fxi}).
\end{remarka}
\begin{remarka} In the Trivialization Lemma \ref{lemma: pi},
we could have worked with the $x-$axis in place of the primary component of the
critical locus. However, when $|a|$ is big, there is no reason to assume that
the $x-$axis will remain transverse to the foliation of $U^+$.
Choosing a transverse which has dynamical
meaning, $\mathcal{C}_0$, gives hope of extending the results to the whole interior
of the hyperbolic component of the \He connectedness locus that
contains perturbations of a hyperbolic polynomial. In fact, Theorem \ref{thm:
LyubichRobertson} is also believed to hold in this generality.
\end{remarka}

\begin{prop}The function $\alpha:\C-\D\rightarrow \C^*$ is well defined, analytic on $\C-\overline{\D}$ and continuous on $\s^1$.
\end{prop}
\proof  For $\xi$ fixed, $\alpha(\xi)$ is defined in Equation \ref{eq:alpha-01} as the
difference quotient of three points $H(c_0(\xi))$, $c_0(\xi^2)$ and
$c_{-1}(\xi^2)$ from $\mathcal{F}_{\xi^2}$. The ratio $(x-y)/(y-z)$ of three
distinct points $x$, $y$, $z$ from a manifold biholomorphic to $\C$ is independent of the
choice a particular trivialization. Hence $\alpha$ is well defined.
The function $c_0:\C-\D\rightarrow \C^2$ is holomorphic on
$\C-\overline{\D}$ and continuous on $\s^1$. The affine structure on
$J^{+}\cup U^{+}$ is transversely continuous \cite{BS5}. 
\qed

\begin{prop}\label{prop:cocycle}
The function $\alpha$ is unique up to multiplication by appropriate
maps of the form $\displaystyle u(\xi^2)/u(\xi)$, where
$u:\C-\D\rightarrow \C^{*}$ is holomorphic on $\C-\overline{\D}$ and
continuous on $\s^1$. In addition, the function
$u(\xi)$ is well defined on the primary component of the critical
locus, that is, if $c_0(\xi_1)=c_0(\xi_2)$ then $u(\xi_1)=u(\xi_2)$.
\end{prop}
\proof Suppose we define another trivialization of $\mathcal{F}_{\xi}$ that
assigns $c_0(\xi)\rightarrow v(\xi)$ and $c_{-1}(\xi)\rightarrow
u(\xi)$, where $u,v:\C-\D\rightarrow \C$ and $u(\xi)\neq v(\xi)$. 
We use two very special transverses to give a trivialization of
$J^+\cup U^+$, namely $\mathcal{C}_0$ and $H^{-1}(\mathcal{C}_0)$,
which carry their own identifications, described in Theorem
\ref{thm: LyubichRobertson}. So the assignments
$c_0(\xi)\rightarrow v(\xi)$ and $c_{-1}(\xi)\rightarrow u(\xi)$
must preserve these identifications, that is, if $c_0(\xi_1)=c_0(\xi_2)$
for some $\xi_1,\xi_2\in \s^1$ then $u(\xi_1)=u(\xi_2)$ and
$v(\xi_1)=v(\xi_2)$. 

 The same computation as before yields
\begin{equation*}
\alpha(\xi)u(\xi)+\beta(\xi)= v(\xi^2)\ \ \ \mbox{and}\ \ \ \alpha(\xi)v(\xi)+\beta(\xi)= x,
\end{equation*}
where $x$ is the coordinate of $H(c_0(\xi))$ in
$\mathcal{F}_{\xi^2}$ and can be computed from the invariance
property of the ratio of three points under affine changes of
coordinates. We have
\[
\frac{H(c_0(\xi))-c_0(\xi^2)}{c_0(\xi^2)-c_{-1}(\xi^2)}=\frac{x-v(\xi^2)}{v(\xi^2)-u(\xi^2)}
\ \Rightarrow\
x=\frac{H(c_0(\xi))-c_0(\xi^2)}{c_0(\xi^2)-c_{-1}(\xi^2)}\left(v(\xi^2)-u(\xi^2)\right)+v(\xi^2).
\]
\noindent After solving the system we get
\begin{eqnarray}\label{eq: 3}
\alpha(\xi) &=&
\frac{v(\xi^2)-x}{u(\xi)-v(\xi)}=\frac{u(\xi^2)-v(\xi^2)}{u(\xi)-v(\xi)}\cdot
\frac{H(c_0(\xi))-c_0(\xi^2)}{c_0(\xi^2)-c_{-1}(\xi^2)}\\
\beta(\xi) &=&
v(\xi^2)-u(\xi)\cdot\frac{u(\xi^2)-v(\xi^2)}{u(\xi)-v(\xi)}\cdot
\frac{H(c_0(\xi))-c_0(\xi^2)}{c_0(\xi^2)-c_{-1}(\xi^2)}
\end{eqnarray}
so the expression of $\alpha$ has only changed by a multiplicative
factor of the form
\begin{eqnarray}\label{eq: uv}
\frac{u(\xi^2)-v(\xi^2)}{u(\xi)-v(\xi)}.
\end{eqnarray}
Since we are looking only for
functions $\alpha(\xi)$ which are analytic on $\C-\overline{\D}$ and
continuous on $\s^1$, the maps $u(\xi)$ and $v(\xi)$ must also be
holomorphic on $\C-\overline{\D}$ and continuous on $\s^1$. We denote the multiplicative
factor \ref{eq: uv} by $u(\xi^2)/u(\xi)$ 
when there is no danger of confusion.
\qed

We will  now show that the Stable Multiplier Condition \ref{eq:restriction} is satisfied.
\begin{prop}\label{prop:product-multipliers}
For $\xi\in \s^1,\ \xi=\xi^{2^k}$, the product
$\alpha(\xi)\alpha(\xi^2)\ldots\alpha(\xi^{2^{k-1}})$
 does not depend on the choices of affine maps and it equals the small
eigenvalue $\lambda(DH_x^{\circ k})$ of the derivative of $H^{\circ
k}$ at some $k$-periodic point $x$ of $H$.
\end{prop}
\proof Let $\xi\in \s^1,\ \xi=\xi^{2^k}$ be a periodic point of the
doubling map $\xi\rightarrow \xi^2$.

The fiber $\mathcal{F}_{\xi}$ is invariant
under the \He map since $H^{\circ
k}\left(\mathcal{F}_{\xi}\right)=\mathcal{F}_{\xi^{2^k}}=\mathcal{F}_{\xi}$,
hence $\mathcal{F}_{\xi}$ is the stable manifold $W^s(x)$ of some
periodic point $x$ of period $k$ of the \He map, 
and $\mathcal{F}_{\xi^{2^i}}$ is the stable
manifold of $H^{\circ i}(x)$, 
$1\leq i\leq k$. Moreover, since $H$ is a hyperbolic \He
map, we know that the tangent space $T_{H^{\circ
i}(x)}\mathcal{F}_{\xi^{2^i}}$ is the eigenspace of the smallest
eigenvalue of the Jacobian matrix $DH^{\circ i}_{H^{\circ i}(x)}$.

The function $\alpha(\xi)$ is unique up to a
multiplicative factor $u(\xi^2)/u(\xi)$. We notice that
\[
\frac{u(\xi^2)}{u^(\xi)}\frac{u(\xi^4)}{u(\xi^2)}\ldots
\frac{u(\xi^{2^{k}})}{u(\xi^{2^{k-1}})}=1,
\]
hence the product
$\alpha(\xi)\alpha(\xi^2)\ldots \alpha(\xi^{2^{k-1}})$ is well
defined and independent of choices of affine maps on the fibers
$\mathcal{F}_{\xi}$, $\mathcal{F}_{\xi^2}$,\ldots,
$\mathcal{F}_{\xi^{2^{k-1}}}$. Therefore
\[
\alpha(\xi)\alpha(\xi^2)\ldots \alpha(\xi^{2^{k-1}})=\prod\limits
_{i=1}^{k}\frac{H(c_0(\xi^{2^{i-1}}))-c_0(\xi^{2^i})}{c_0(\xi^{2^i})-H^{-1}\left(
c_{0}(\xi^{2^{i+1}})\right)},
\]
where each of the ratios is evaluated in
$\mathcal{F}_{\xi^{2^{i}}},1\leq i\leq k$.

\noindent Each fiber is biholomorphic to $\C$ and we
can choose convenient parametrizing functions
\[
\psi_i:\C\rightarrow \mathcal{F}_{\xi^{2^i}}\ \ \mbox{with}\ \
\psi_i(0)=H^{\circ i}(x)\ \ \mbox{and}\ \ \psi_i'(0)=v_i,
\]
where $v_0$ is a stable eigenvector of $DH_x$ and $v_i=DH^{\circ
i}_{H^{\circ i}(x)}v$.
 Denote by
$\phi_i:\mathcal{F}_{\xi^{2^i}}\rightarrow \C$ the inverse function
of $\psi_i$.
\[
\diag{\mathcal{F}_{\xi}}{\mathcal{F}_{\xi^2}}{\C}{\C}{H}{\phi_0}{\phi_1}{L_1(z)=m_1\cdot
z}\ \ \ \ \ldots \ \ \ \
\diag{\mathcal{F}_{\xi_{2^{k-1}}}}{\mathcal{F}_{\xi^{2^k}}}{\C}{\C}{H}{\phi_{k-1}}{\phi_k=\phi_0}{L_k(z)=m_k\cdot
z}
\]

\vglue 0.25cm
\noindent The \He map induces multiplicative maps
between the copies of $\C$, $\phi_i\circ H = m_i\cdot \phi_{i-1}$
where $m_i\neq 0$ and we have
\[
\phi_0\circ H^{\circ k}= \left(m_k\cdot m_{k-1}\cdot \ldots \cdot
m_1\right) \cdot \phi_0.
\]
\noindent If we differentiate the previous relation and evaluate at
$x$ we get
\[
\triangledown \phi_0\cdot DH^{\circ k}_{x}\cdot v_0= \left(m_k\cdot
m_{k-1}\cdot \ldots \cdot m_1\right) \cdot \triangledown \phi_0\cdot
v_0.
\]
\noindent But $DH^{\circ k}_{x}\cdot v_0=\lambda v_0$, where
$\lambda$ with $|\lambda|<1$ is the small eigenvalue of the Jacobian
matrix $DH^{\circ k}_{x}$. Hence
\[
m_k\cdot m_{k-1}\cdot \ldots \cdot m_1=\lambda.
\]
\noindent One can now compute the product
\[\prod\limits_{i=1}^k \frac{\phi_i\circ H(c_0(\xi^{2^{i-1}}))-\phi_i\circ c_0(\xi^{2^i})}{\phi_i\circ
c_0(\xi^{2^i})-\phi_i\circ
H^{-1}\left(c_{0}(\xi^{2^{i+1}})\right)}=\prod\limits_{i=1}^k
\frac{\phi_i\circ H(c_0(\xi^{2^{i-1}}))-\phi_i\circ H \circ H^{-1}
(c_0(\xi^{2^i}))}{\phi_i\circ c_0(\xi^{2^i})-\phi_i\circ
H^{-1}\left(c_{0}(\xi^{2^{i+1}})\right)}\]
\begin{eqnarray*}
&=& \prod\limits_{i=1}^k \frac{m_i\cdot \phi_{i-1}\circ
c_0(\xi^{2^{i-1}})-m_i\cdot \phi_{i-1}\circ H^{-1}(c_0(\xi^{2^i}))}{
\phi_i\circ c_0(\xi^{2^i})-\phi_i\circ
H^{-1}\left(c_{0}(\xi^{2^{i+1}})\right) }\\
&=&(m_1\cdot m_2\cdot\ldots \cdot m_k) \frac{\phi_0\circ
c_0(\xi)-\phi_0\circ H^{-1}(c_0(\xi^2))}{\phi_k\circ
c_0(\xi^{2^k})-\phi_k\circ H^{-1}(c_0(\xi^{2^{k+1}}))}=\lambda.
\end{eqnarray*}
\vglue -0.5cm
\qed

The description of candidate functions $\alpha(\xi)$ from
proposition \ref{prop:cocycle} that satisfy the condition in
\ref{prop:product-multipliers} can be linked to other results of
this sort.

\begin{thm}[Livschitz \cite{K}]\label{thm:Livschitz}
If $\Lambda$  is a topologically transitive hyperbolic set for a
diffeomorphism $f$  and $\varphi:\Lambda\rightarrow \R$  is a
$\tau$-H\"{o}lder continuous function such that
\[
\sum^{n-1}_{i=0} \varphi(f^i (x))=0 \mbox{ whenever } f^n(x)=x,
\]
then $\varphi$ is a coboundary, i.e. there exists a continuous
function $h:\Lambda\rightarrow\R$ such that $\varphi=h\circ f-h$.
This function is unique up to an additive constant, and it is a
$\tau$-H\"{o}lder continuous.
\end{thm}
 \begin{prop}\label{prop:Lamin1}
Suppose $\xi_1, \xi_2$ are two points on $\s^1$ such that
$\gamma(\xi_1)=\gamma(\xi_2)$, where $\gamma$ is the Carath\'eodory
loop of the polynomial $p$. Then $\alpha(\xi_1)=\alpha(\xi_2)$.
\end{prop}
\proof We first show that
\begin{eqnarray}\label{eq: 2}
\frac{H(c_0(\xi_1))-c_0(\xi_1^2)}{c_0(\xi^2_1)-H^{-1}(c_0(\xi^4_1))}
=\frac{H(c_0(\xi_2))-c_0(\xi_2^2)}{c_0(\xi^2_2)-H^{-1}(c_0(\xi^4_2))}.
\end{eqnarray}

Since $c_0(\xi)= \left(\tau^{+}|_{\overline{\mathcal{C}}_0}\right
)^{-1}\circ \gamma(\xi)$ and  $\gamma(\xi_1)=\gamma(\xi_2)$ we have
$c_0(\xi_1)=c_0(\xi_2)$. By using the properties of the B\"{o}ttcher
coordinate $\gamma(\xi^2)=p(\gamma(\xi))$, $\xi\in\C-\D$, we
get that $c_0(\xi_1^2)=c_0(\xi_2^2)$ and $c_0(\xi_1^4)=c_0(\xi_2^4)$.
Hence the two ratios in Equation \ref{eq: 2} are equal.

By equation \ref{eq: 3} in Proposition \ref{prop:cocycle}, the
choice of $\alpha$ is unique up to multiplication by functions of
the form $u(\xi^2)/u(\xi)$. These functions $u$ satisfy the additional
property that $c_0(\xi_1)=c_0(\xi_2)\Rightarrow u(\xi_1)=u(\xi_2)$.
By Theorem \ref{thm: LyubichRobertson}, $c_0(\xi_1)=c_0(\xi_2)$ if and only
if $\gamma(\xi_1)=\gamma(\xi_2)$.
\qed

\begin{prop}\label{prop: ergodic} $\int_{\s^1}\log|\alpha(e^{2 \pi i \theta})|d\theta = \lambda^{-}(\mu)$, 
where
$\lambda^{-}(\mu)$ is the stable Lyapunov exponent with respect to
the unique measure $\mu$ of maximal entropy supported on the Julia set $J$.
\end{prop}
\proof The stable and unstable Lyapunov exponent are well understood
in the case of hyperbolic \He maps. They are related by the equation
\[
\lambda^{+}(\mu)+\lambda^{-}(\mu)=\log|DH|=\log|a|.
\]
When $|a|<1$ and the \He map is hyperbolic with connected Julia set,
the unstable Lyapunov exponent is $\lambda^{+}(\mu)=\log(2)$, as
shown in \cite{BS5}. Hence $\lambda^{-}(\mu)=\log|a|-\log 2$. The
stable Lyapunov exponent $\lambda^{-}$ is defined as
\[
\lambda^{-}=\lim\limits_{k\rightarrow \infty} \frac{1}{k}\log
\|DH^{\circ k}|_{E^s}\|.
\]
By \cite{BS5}, for $\mu$ almost every point $x$ in $J$, one has
\[
\lambda^{-}=\lim\limits_{k\rightarrow \infty} \frac{1}{k}\log
\|DH^{\circ k}_x|_{E_x^s}\|.
\]
Let $\mathcal{F}_{\xi}$ be a leaf of the lamination of $J^+$ that
contains $x$. We  can compute $\lambda^{-}$ as follows
\begin{eqnarray*}
\lambda^{-}&=&\lim\limits_{k\rightarrow \infty} \frac{1}{k}\log
\|DH^{\circ k}_x|_{E_x^s}\|=\lim\limits_{k\rightarrow \infty}
\frac{1}{k}\log\big{|}\alpha(\xi)\alpha(\xi^2)\cdot\ldots \cdot
\alpha(\xi^{2^{k-1}})\big{|} \\
&=& \lim\limits_{k\rightarrow \infty}
\frac{1}{k}\left(\log|\alpha(\xi)|+\log|\alpha(\xi^2)|+\ldots +
\log|\alpha(\xi^{2^{k-1}})|\right)=\int_{\s^1}\log|\alpha(\xi)|d\theta.
\end{eqnarray*}
Here $\int_{\s^1}\log|\alpha(\xi)|d\theta$ stands
for $\int_{\s^1}\log|\alpha(e^{2 \pi i \theta})|d\theta$,
where $d\theta$ is the Lebesgue measure on the unit circle $\s^1$ regarded here as $\R/\Z$.
The doubling map $f(\xi)=\xi^2$ is ergodic with respect to the
Lebesgue measure on $\s^1$, so the orbit of almost every $\xi$ is
equidistributed on $\s^1$. The last equality then follows from the
Birkhoff Ergodic Theorem.
\qed
\begin{remarka}
Notice that $\int_{\s^1}\log|\alpha(\xi)|d\theta$ in Lemma \ref{prop: ergodic} does not depend on the choices involved in the construction
of the function $\alpha$. The map $\alpha$ is unique up to
multiplication by a factor of the form $u(\xi^2)/u(\xi)$, where the map
$u:\s^1\rightarrow \C^*$ is continuous. Since $f(\xi)=\xi^2$ is
measure preserving with respect to the Lebesgue measure on $\s^1$,
we have
$\int_{\s^1}\log|u(\xi^2)|d\theta=\int_{\s^1}\log|u(\xi)|d\theta$.
\end{remarka}

The function $\alpha$ is probably a full invariant of the
(quadratic) \He map, in the sense that if two hyperbolic  \He maps
$H_1$ and $H_2$ have the property that $\alpha_1=\alpha_2$ then
the \He maps coincide, i.e. $H_1=H_2$. 
The following proposition from \cite{T} provides support for this
claim.

\begin{prop} The values $\alpha(1)$ and $\alpha(e^{2\pi i 1/3})\cdot\alpha(e^{2\pi i 2/3})$
determine the \He map up to three choices.
\end{prop}
Moreover, by Proposition \ref{prop: ergodic} we have $\log|a|-\log(2)=\int_{\s^1}\log|\alpha(e^{2 \pi i \theta})|d\theta$, therefore the absolute value of the Jacobian is determined by the function $\alpha$. It would also be interesting to study the relation between the function $\alpha$ and the non-transversality locus invariant (ntl-invariant) described in \cite{HOV3}.

\section{Degeneracy of the function $\alpha$}\label{Degeneracy}

It is easy to see that the limit of the function $\alpha$ is zero when
the Jacobian goes to zero. This is a consequence of the
fact that the critical points on the primary component
$\mathcal{C}_0$ and on $\mathcal{C}_{1}=H(\mathcal{C}_0)$ remain
bounded as $a\rightarrow 0$ and close to the $x$-axis,
respectively to the parabola $\{(x,y)\in \C^2, x=p(y)\}$. Meanwhile, by the definition of $H^{-1}$ from Section \ref{sec:Preliminaries}, 
the critical points on $\mathcal{C}_{-1}=H^{-1}(\mathcal{C}_0)$ go
to infinity as the Jacobian tends to $0$.

It is therefore more interesting and useful
to compute the limit of $\alpha(\xi)/a$ as $a\rightarrow 0$.
 In the trivialization that assigns $c_0(\xi)\rightarrow
0$ and $c_{-1}(\xi)\rightarrow 1$ we have computed in
Equation \ref{eq:alpha-01} the following formula 
\begin{equation}\label{eq:alpha-loc}
    \alpha(\xi)=\frac{H(c_0(\xi))-c_0(\xi^2)}{c_0(\xi^2)-H^{-1}(c_0(\xi^4))}.
\end{equation}
Throughout this section, we will refer to trivialization \ref{eq:sttriv} as the
{\it standard trivialization}  and to the function $\alpha$ from  Equation \ref{eq:alpha-loc} as
the function $\alpha$ computed with respect to the standard
trivialization.

\begin{prop}[\textbf{Degeneracy of $\alpha$}]\label{Prop:triv1}
\[
\displaystyle{\lim\limits_{a\rightarrow
0}\frac{1}{a}\cdot\frac{H(c_0(\xi))-c_0(\xi^2)}{c_0(\xi^2)-H^{-1}(c_0(\xi^4))}}=
\frac{\gamma(\xi)}{2(\gamma(\xi^2))^2}
\]
where $\gamma$ is the B\"ottcher  coordinate of the polynomial $p$.
\end{prop}
\proof
Let $\xi\in \C-\D$ and set $x=c_0(\xi)$. The leaf
$\mathcal{F}_{\xi}$ is isomorphic to $\C$ and one can choose a
biholomorphic map
\[
\psi_x:\C\rightarrow
\mathcal{F}_{\xi}
\]
such that $\psi_x(0)=x$ and $\psi_x'(0)=v$, where $v\in
T_x\mathcal{F}_{\xi}$ is the unit vector $v$, with $\|v\|=1$ and
$pr_2(v)\in \R^+$.

We know that when $a$ is small, $\mathcal{F}_{\xi}$ is almost
vertical in a neighborhood of $x$ \cite{HOV2} \cite{T}. Hence there exists a
unique analytic function $f:\D\rightarrow \C$ such that locally
around the critical point $x=(f(z_0),z_0)$, the leaf $\mathcal{F}_{\xi}$ is the graph of $f$, of the form 
$\{(f(z),z)\}$. So we can choose $v=(f'(z_0),1)/\|(f'(z_0),1)\|$.

The leaves $\mathcal{F}_{\xi^2}$ and
$\mathcal{F}_{\xi^{2^2}}$ are isomorphic to $\C$, hence there exist biholomorphisms
\[
\psi_{i}:\C\rightarrow \mathcal{F}_{\xi^{2^i}},\ \ i=1,2
\]
\noindent such that $ \psi_{i}(0)=H^{\circ i}(x)\ \mbox{and}\
\psi_{i}'(0)=w_i,\ \mbox{where}\ w_i\in T_{H^{\circ
i}(x)}\mathcal{F}_{\xi^{2^i}}$ is a tangent vector with norm
$\|w_i\|=1$ and $pr_2(w_i)\in \R^+$.
\[
\minCDarrowwidth60pt
\begin{CD}
{\mathcal{F}_{\xi}} @> {H} >> {\mathcal{F}_{\xi^2}} @>{H}>> {\mathcal{F}_{\xi^4}}\\
@A{\psi_x}AA  @AA{\psi_{1}}A @AA{\psi_{2}}A\\
{\C} @>{L_1(z)=m_1\cdot z}>> {\C} @>{L_2(z)=m_{2}\cdot z}>> {\C}
\end{CD}
\]

\noindent From this commutative diagram we get
\begin{eqnarray*}
H\circ \psi_{1}(z)&=&\psi_{2}(m_{2}\cdot z)\\
\psi^{-1}_{1}\circ H^{-1}&=&\frac{1}{m_{2}}\cdot \psi_{2}^{-1}.
\end{eqnarray*}
We can therefore compute the function $\alpha(\xi)$ using Equation \ref{eq:alpha-loc} and get
\[
\alpha(\xi)=\frac{\psi^{-1}_{1}\circ H(c_0(\xi))-\psi^{-1}_{1}\circ
c_0(\xi^2)}{\psi^{-1}_{1}\circ c_0(\xi^2)-\psi^{-1}_{1}\circ
H^{-1}(c_0(\xi^4))}=\frac{\psi^{-1}_{1}\circ
H(c_0(\xi))-\psi^{-1}_{1}\circ c_0(\xi^2)}{\psi^{-1}_{2}\circ H(
c_0(\xi^2))-\psi^{-1}_{2}\circ c_0(\xi^4)}\cdot m_{2}.
\]
\noindent When the Jacobian $a$ is $0$, the primary component of the
critical locus degenerates uniformly to $(\C-\mathring{K}_p)\times \{0\}$,
where $K_p$ is the filled-in Julia set of the polynomial
$p$. Morover $c_0(\xi)=(\gamma(\xi),0)$,
$c_0(\xi^2)=(\gamma(\xi^2),0)$ and $c_0(\xi^4)=(\gamma(\xi^4),0)$.
With our notation $x=c_0(\xi)$, we can also compute the degeneracy
of the points $H(x)$ and $H^{\circ 2}(x)$:
\begin{eqnarray*}
H(x) &=&
(\gamma(\xi)^2+c,\gamma(\xi))=(\gamma(\xi^2),\gamma(\xi))\\
H^{\circ 2}(x) &=&
(\gamma(\xi^2)^2+c,\gamma(\xi^2))=(\gamma(\xi^4),\gamma(\xi^2)).
\end{eqnarray*}
\noindent The leaf $\mathcal{F}_{\xi^2}$ degenerates to a collection
of vertical lines
\[
\left\{c_0\left(\omega \xi^2\right)\times \C\ : \ \omega^{2^i}=1, \mbox{ for
some integer } i\geq 1\right\}\!.
\]
However, the parametrizing
function $\psi_{1}$ degenerates to the parametrization of the
vertical line that passes through $c_0(\xi^2)$,
\[
\psi_{1}:\C\rightarrow \gamma(\xi)\times \C,\ \
\psi_{1}(z)=\left(\gamma(\xi^2),z+\gamma(\xi)\right)\!.
\]
\noindent The parametrizing function $\psi_{2}$ degenerates to
\[
\psi_{2}:\C\rightarrow \gamma(\xi^2)\times \C,\ \
\psi_{2}(z)=\left(\gamma(\xi^4),z+\gamma(\xi^2)\right)\!.
\]
\noindent Hence
\begin{eqnarray}\label{m0}
\lim\limits_{a\rightarrow 0}\frac{\psi^{-1}_{1}\circ
H(c_0(\xi))-\psi^{-1}_{1}\circ c_0(\xi^2)}{\psi^{-1}_{2}\circ H(
c_0(\xi^2))-\psi^{-1}_{2}\circ
c_0(\xi^4)}=\frac{\gamma(\xi)}{\gamma(\xi^2)}.
\end{eqnarray}
\noindent From the commutative diagram we also know that
\[
DH_{\psi_1(z)}\cdot \psi_1'(z)=m_2\cdot \psi_{2}'(m_2\cdot z).
\]
When $z=0$ we have $DH_{H(x)}\cdot \psi_1'(0)=m_2 \cdot
\psi_{2}'(0)$, or equivalently
\begin{eqnarray}\label {m2}
\psi_1'(0)=m_2 \cdot DH_{H(H(x))}^{-1}\cdot \psi_{2}'(0).
\end{eqnarray}
\noindent The \He map is $H(x_1,x_2)=\left(x_1^2+c-ax_2,x_1\right)$ and the
inverse has the formula
$H^{-1}(x_1,x_2)=\left(x_2,(x_2^2+c-x_1)/a\right)$, so
\[
DH_{x}^{-1}=\left(
                 \begin{array}{cc}
                   0 & 1 \\
                   \small{-1/a} & \small{2x_2/a} \\
                 \end{array}
               \right).
\]
\noindent It follows that
\[
DH_{H(H(x))}^{-1}=\left(
                 \begin{array}{cc}
                   0 & 1 \\
                   -\small{1/a} & \small{2(x_1^2+c-ax_2)/a} \\
                 \end{array}
               \right)
\]
and Equation \ref{m2} becomes
\begin{eqnarray}\label{m3}
\psi_1'(0)=\frac{m_2}{a} \cdot \left(
                 \begin{array}{cc}
                   0 & a \\
                   -1 & 2(x_1^2+c-ax_2) \\
                 \end{array}
               \right) \cdot \psi_{2}'(0).
\end{eqnarray}
 Notice also that $\lim\limits_{a\rightarrow 0}\psi_1'(0)=\lim\limits_{a\rightarrow 0}\psi_2'(0)=(0,1)$. Thus from Equation \ref{m3} we get
\begin{eqnarray}\label{m4}
\lim\limits_{a\rightarrow 0}\frac{m_2}{a}=\frac{1}{2\gamma(\xi^2)}.
\end{eqnarray}
 Therefore from the relations \ref{m4} and \ref{m0} we can conclude
that
\[
\lim\limits_{a\rightarrow
0}\frac{\alpha(\xi)}{a}=
\frac{\gamma(\xi)}{\gamma(\xi^2)}\cdot\frac{1}{2\gamma(\xi^2)}=\frac{\gamma(\xi)}{2(\gamma(\xi^2))^2}.
\]
 An important observation is that the convergence is
uniform in $\xi$. This follows as a consequence of the fact that the
primary component of the critical locus moves holomorphically with
respect to $a$ when $a$ is small \cite{LR} and degenerates uniformly when $a$
goes to $0$ to $(\C-\mathring{K_p})\times \C$. 
\qed

A consequence of Lemma \ref{Degeneracy} is that the
argument of the function $\alpha|_{\s^1}(\xi)$, regarded as a
function from $\s^1$ to $\s^1$, has degree $-3$. This makes the
plots of the image of the function $\alpha|_{\s^1}$ hard to read.
The following lemma provides a remedy.

\begin{prop}\label{Prop:triv2}
One can choose an appropriate trivialization so that
\[
\lim\limits_{a\rightarrow 0}\frac{\alpha(\xi)}{a}=\frac{1}{2\gamma(\xi)}.
\]
\end{prop}
\proof
Define a trivialization of $\mathcal{F}_{\xi}$ that assigns
$c_0(\xi)\rightarrow 0$ and $c_{-1}(\xi)\rightarrow
\gamma^2(\xi)$, where $\gamma$ is the B\"otcher isomorphism of $p$. Note
that this is an allowed assignment since it verifies the restrictions
in Proposition \ref{prop:cocycle}.
\qed

The insight of Lemma \ref{Prop:triv2} is that $\alpha$ measures the
contraction induced by the derivative of the \He map on the leaves
of the lamination of $J^+\cup U^+$, whereas
\[
2\gamma(\xi)=p'(\gamma(\xi))
\]
measures the expansion of the polynomial $p$ on the Julia set $J_p$
and on $\C-K_p$. As the Jacobian $a$ becomes small, these two
quantities behave like the small and respectively the big eigenvalue
of the Jacobian $DH$ of the \He map, so their product is close to
the determinant $\mbox{det}(DH)=a$.

\section{The image of the cocyle $\alpha$ on the unit circle}\label{Algorithm}

The most interesting behavior of the function $\alpha$ is on the
unit circle $\s^1$. We know that $\alpha:\C-\D\rightarrow \C^{*}$ is
continuous on $\s^1$, nonetheless, we expect that $\alpha$ gives rise to a fractal set when restricted to $\s^1$.  
In \cite{T}, we have designed an algorithm in Python for computing the image of $\alpha$ on the unit circle. Here are some pictures obtained with our program.

\begin{figure}[htb]
\begin{center}
\mbox{\subfigure{ 
\includegraphics[scale =0.48, bb= 190 0 600 375]{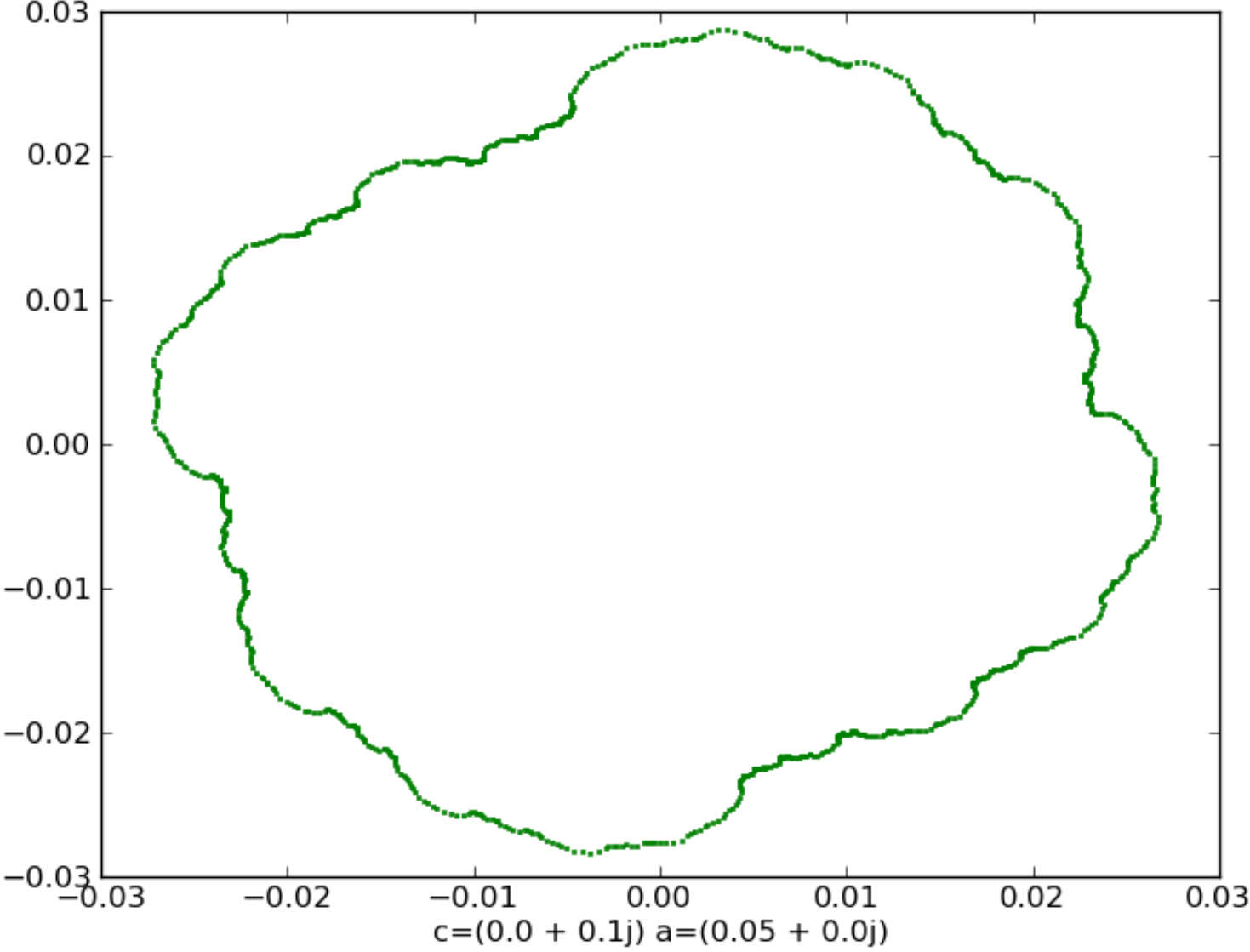} 
\quad 
\hspace{-0.5cm}
\includegraphics[scale=0.41,bb= 225 16 300 430]{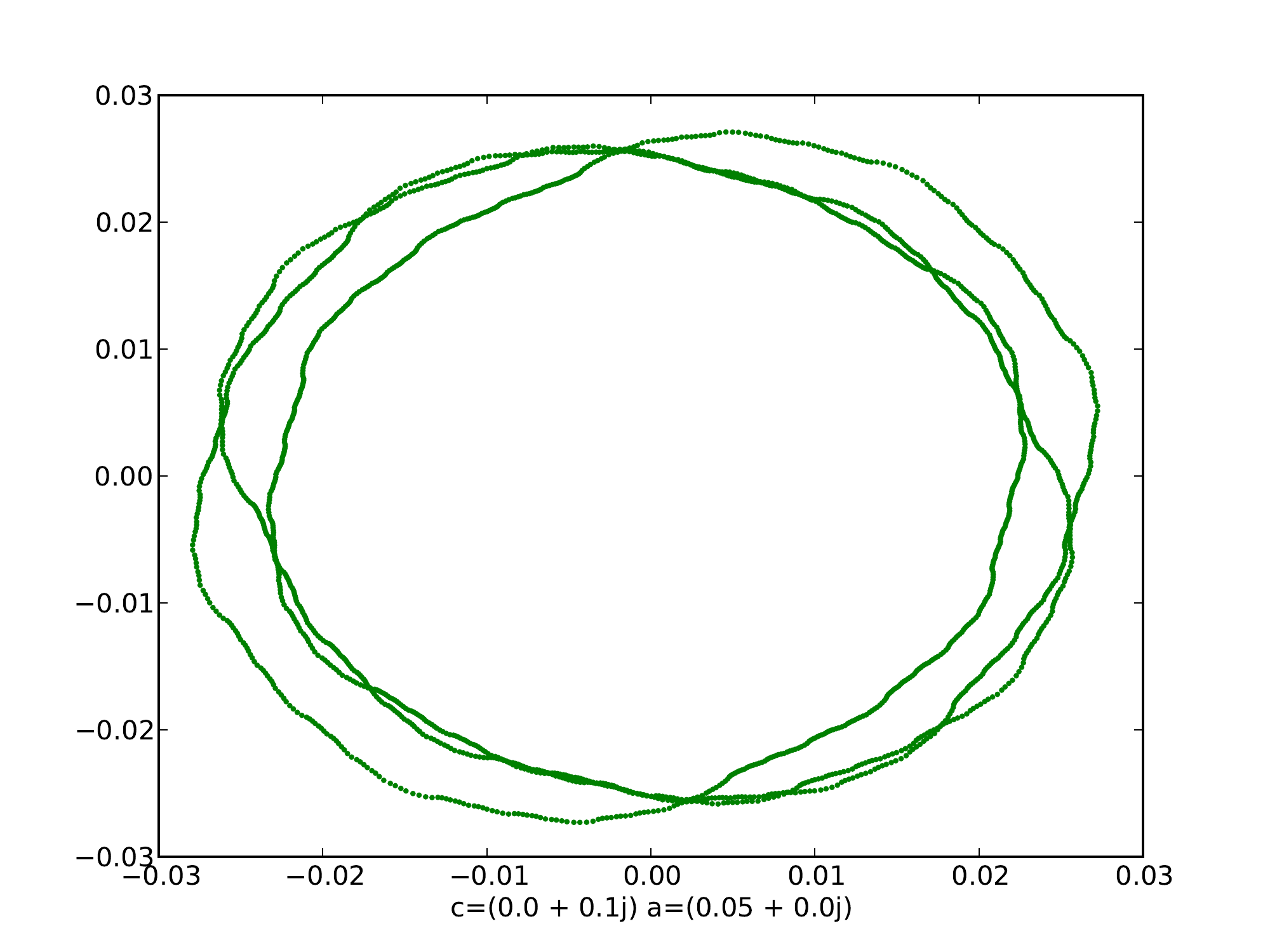}
}
}
\mbox{\subfigure{
\includegraphics[scale=0.403, bb= 255 0 600 430]{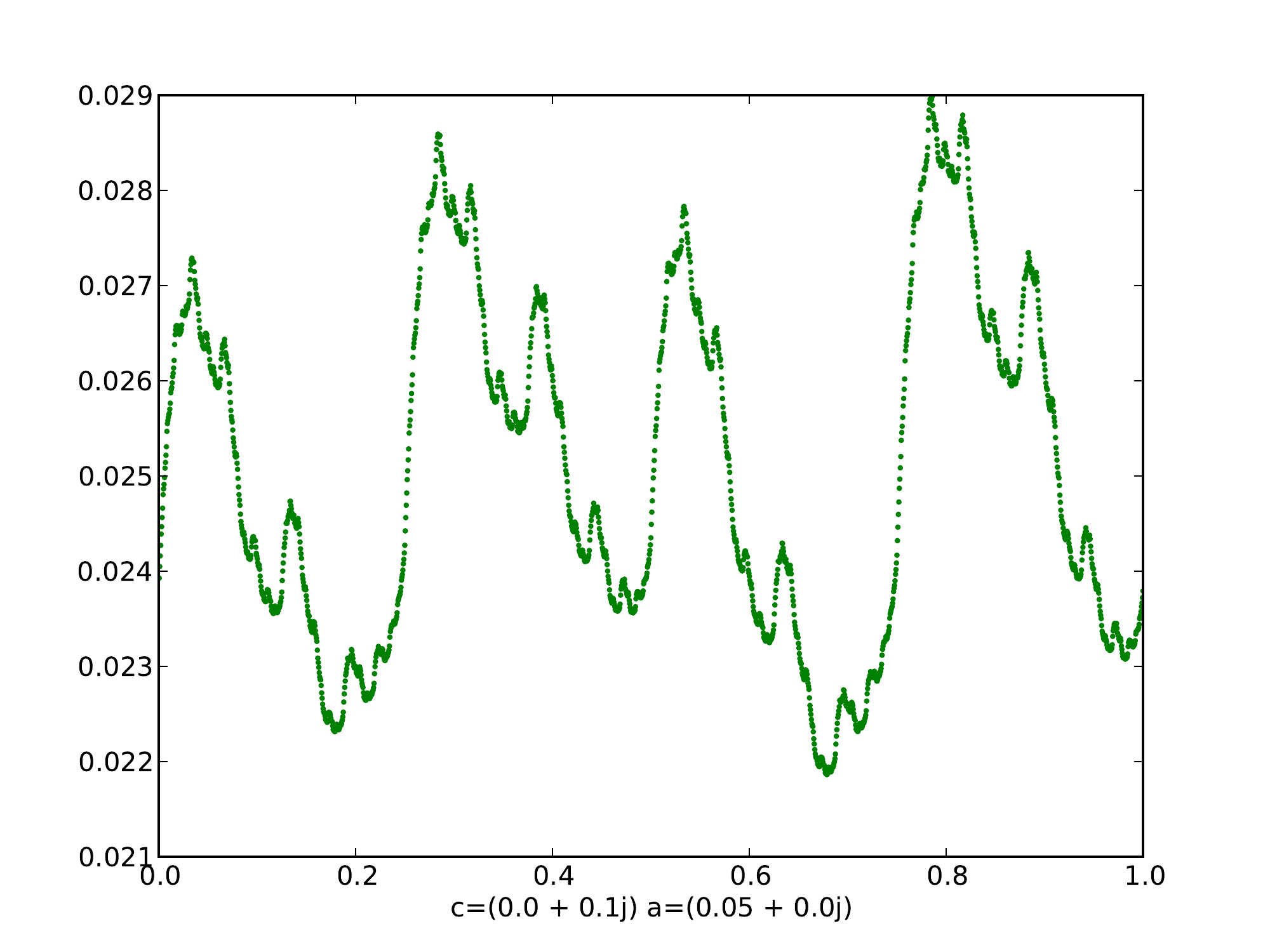}
}
\quad 
\hspace{-0.5cm}
\subfigure{
\includegraphics[scale=0.403, bb= 100 0 300 430]{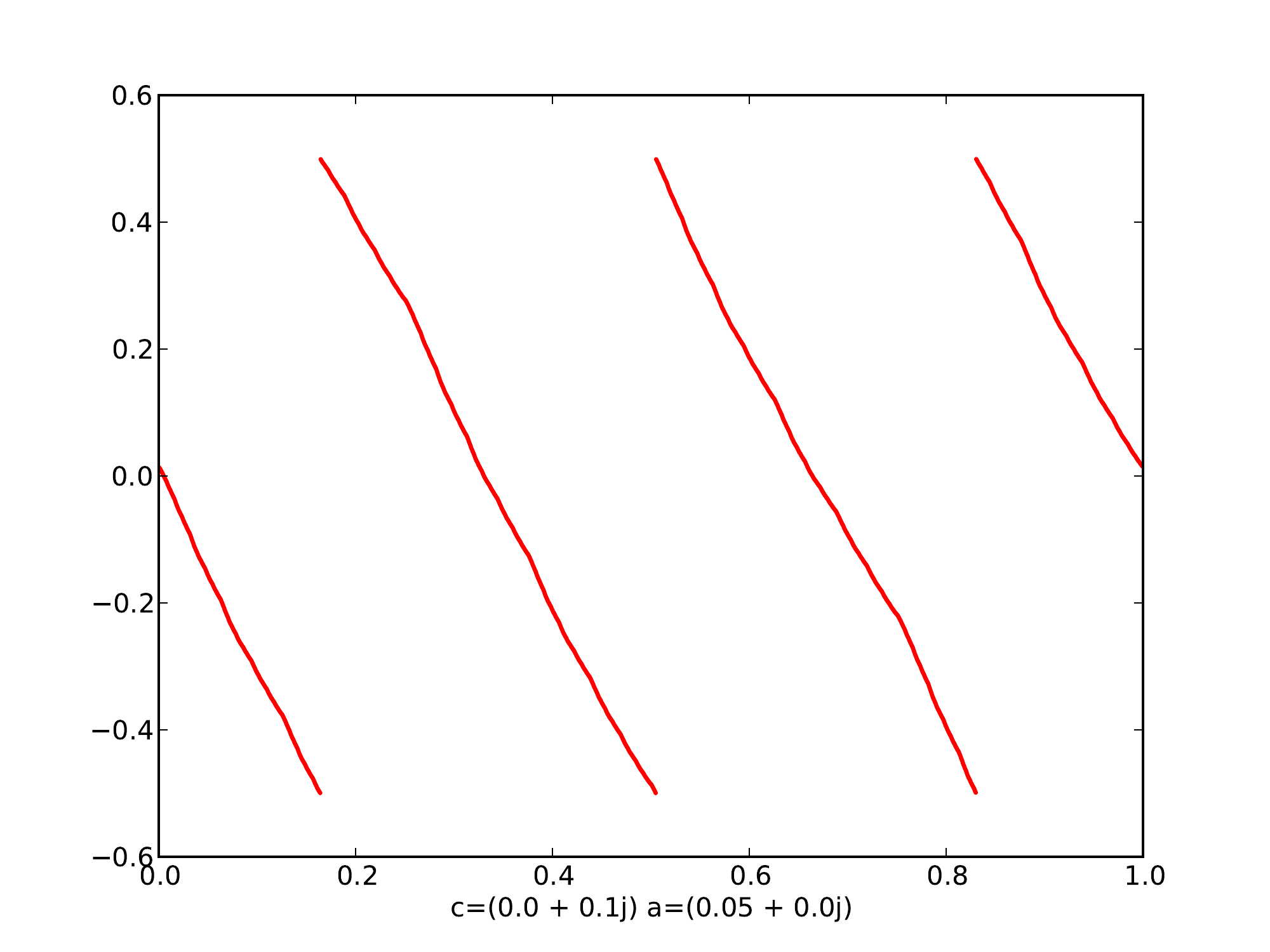} 
}}
\end{center}
\caption{Pictures for the parameters $c=0.1i$ and $a=0.05$. {\sc Top Left:} The image of $\alpha$ on $\s^1$ ($\alpha$ is computed with respect to the trivialization from Proposition \ref{Prop:triv2}). {\sc Top Right:} The image of $\alpha$ on $\s^1$ ($\alpha$ is computed with respect to the standard trivialization). {\sc Bottom Left:}
The graph of the absolute value of $\alpha$ on $\s^1$ (standard trivialization) . {\sc Bottom Right:}
The graph of the argument of $\alpha$ on $\s^1$ (standard trivialization). The argument is
regarded as a function from $[0,1]$ to $[-1/2, 1/2]$.}
\label{fig:im-11-ModArg}
\end{figure}

\begin{figure}[htb]
\begin{center}
\subfigure{\includegraphics[scale=0.57, bb= 290 -400 290 430]{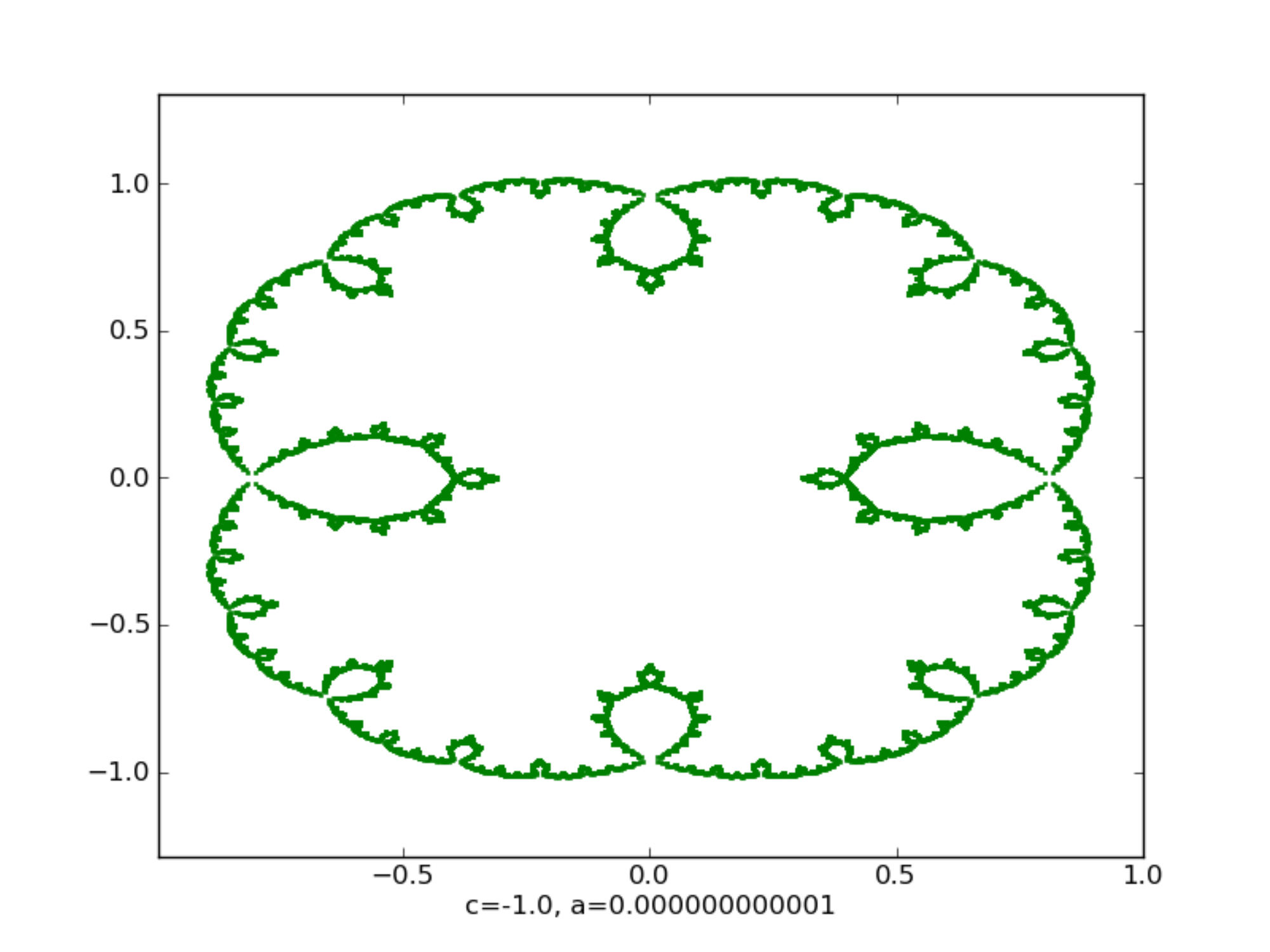} }
\subfigure{\includegraphics[scale=0.57, bb= 300 20 300 430]{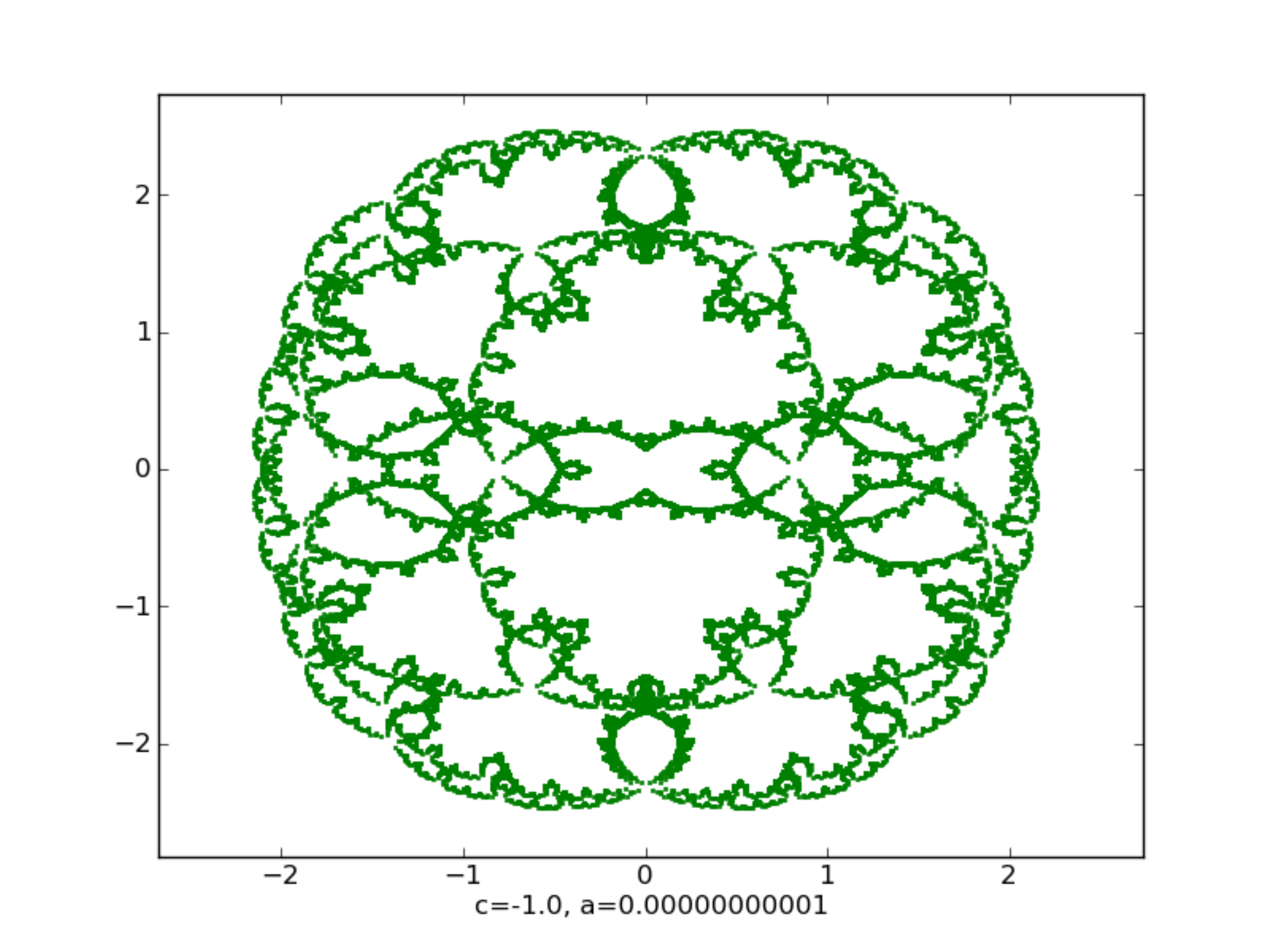}  }
\end{center}
\caption{Here $c=-1$ and $a$ is very small. {\sc Top:} The image of
$\alpha/a$ on $\s^1$ with the trivialization from Proposition \ref{Prop:triv2}.
{\sc Bottom:} The image of $\alpha/a$ on $\s^1$ with the standard
trivialization from Proposition \ref{Prop:triv2}. Regarded as a function on $\s^{1}$, $\alpha$ has
degree $-3$. } 
\label{fig:im-20}
\end{figure}

\begin{figure}[htb]
\begin{center}
\mbox{
\subfigure{\includegraphics[scale=0.39, bb= 500 0 300 430]{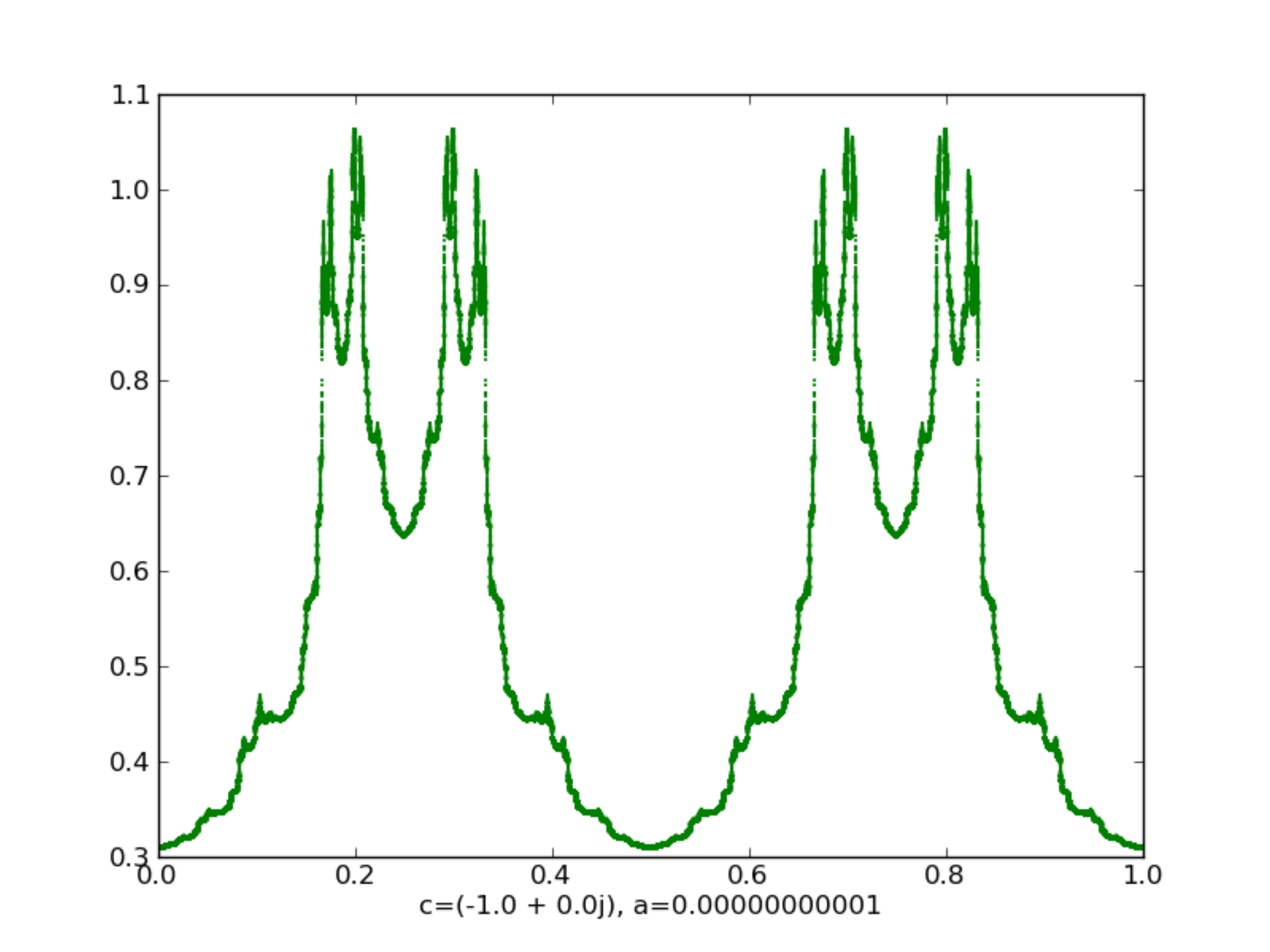}} 
\quad
\subfigure{\includegraphics[scale=0.39, bb= 0 0 70 430]{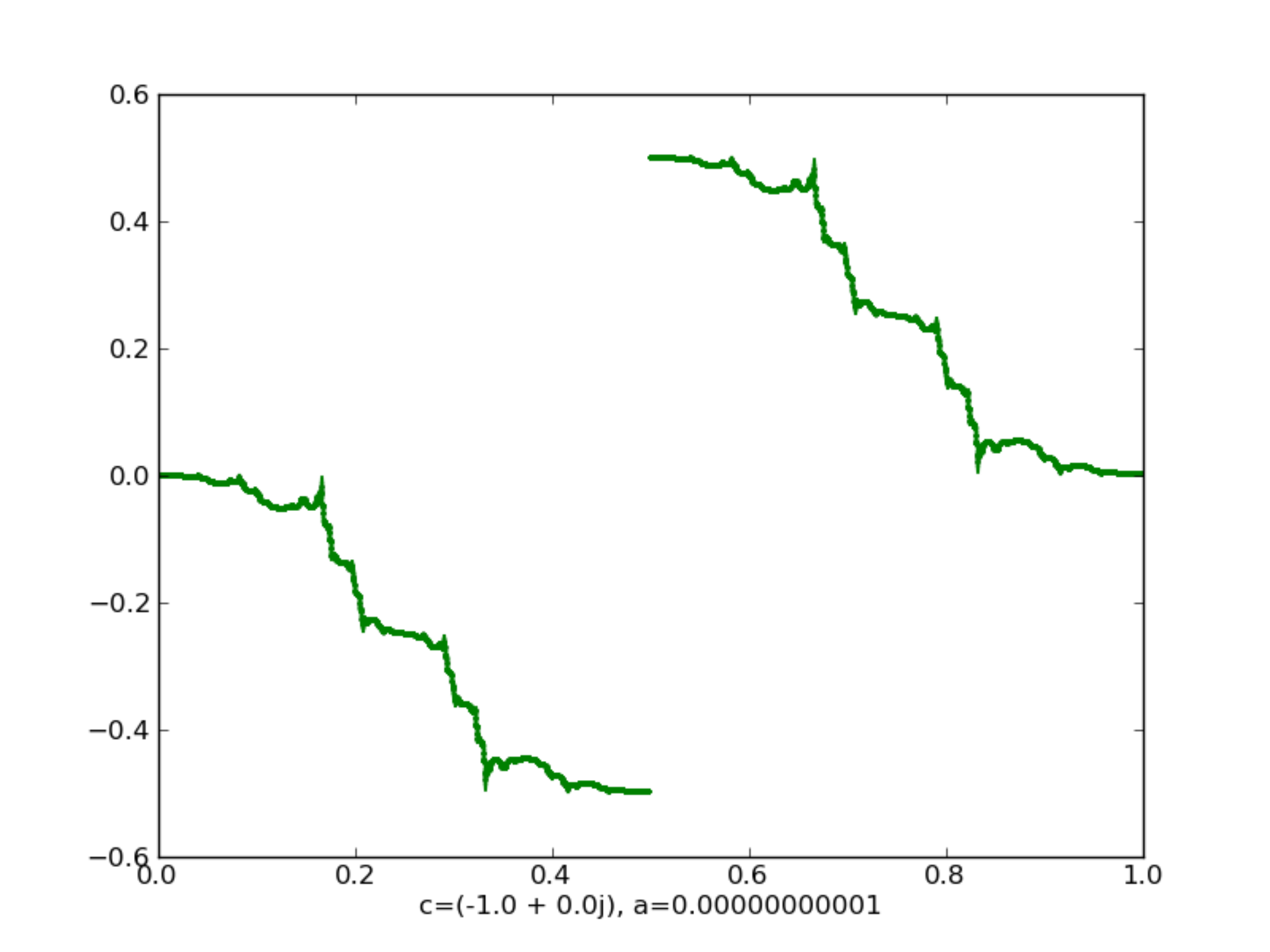} }}
\end{center}
\caption{Pictures for the parameters $c=-1$ and $a=10^{-10}$. The
function $\alpha$ has been computed with respect to the
trivialization from \ref{Prop:triv2}.  {\sc Left:} The graph of the
absolute value of $\alpha/a$ on $\s^1$. {\sc Right:} The graph of the
argument of $\alpha/a$ on $\s^1$. The argument is regarded as a
function from $[0,1]$ to $[-1/2,1/2]$.}
\label{fig:im-19-ModArg}
\end{figure}

\section{Growth estimates for the group $\Gamma_{p,a}$}\label{sec:grgamma}
The main motivation for studying the properties of $\alpha$ on $\s^1$ is that all other extensions can be expressed as functions of $\alpha$. The lift $\widetilde{H}$ of the \He map, as well as the elements of the group of deck transforms  $\Gamma_{p,a}$ can be extended to $\s^1\times\C$ and we can now recursively compute the elements of the group 
\begin{eqnarray}\label{group}
\Gamma_{p,a}=\left\{
\gamma_{\frac{j}{2^{k}}} \pvec{\xi}{z}=\pvec{e^{2\pi i\frac{j}{2^k}}\xi}{p_{j,k}(\xi)z+q_{j,k}(\xi)}, 1\leq j\leq
2^k\right\}
\end{eqnarray}
from the relation
$\widetilde{H}\circ\gamma_{\frac{j}{2^{k+1}}}=\gamma_{\frac{j}{2^{k}}}\circ\widetilde{H}$
and $\gamma_1=\mbox{id}$ and get
\begin{eqnarray*}\label{eq:Group Recursive}
p_{j,k}(\xi)&=& p_{j,k-1}(\xi^2)\frac{\alpha(\xi)}{\alpha\left(e^{2\pi i\frac{j}{2^k}}\xi\right)}\\
q_{j,k}(\xi)&=&\frac{p_{j,k-1}(\xi^2)\beta(\xi)+q_{j,k-1}(\xi^2)-\beta\left(e^{2\pi
i\frac{j}{2^k}}\xi\right)} {\alpha\left(e^{2\pi
i\frac{j}{2^k}}\xi\right)}.
\end{eqnarray*}
Indeed, we know that $\widetilde{H}(\xi,z)=(\xi^2,
\alpha(\xi)z+\beta(\xi))$, $\gamma_1(\xi,z)=(\xi,z)$, and for $k>1$
\[
\widetilde{H}\circ \pvec{e^{2\pi
i\frac{j}{2^k}}\xi}{p_{j,k}(\xi)z+q_{j,k}(\xi)}=\pvec{e^{2\pi
i\frac{j}{2^{k-1}}\xi}} {p_{j,k-1}(\xi)z+q_{j,k-1}(\xi)}\circ
\widetilde{H}\hvec{\xi}{z}.
\]
By comparing the second coordinate we get the following relation
\[
\alpha\left(e^{2\pi i\frac{j}{2^k}}\xi\right)\cdot
\left(p_{j,k}(\xi)z+q_{j,k}(\xi)\right)+\beta(e^{2\pi
i\frac{j}{2^k}}\xi)=p_{j,k-1}(\xi^2)\left(\alpha(\xi)z+\beta(\xi)\right)+q_{j,k-1}(\xi^2).
\]
Therefore
\[
\left\{
  \begin{array}{ll}
    \alpha\left(e^{2\pi i\frac{j}{2^k}}\xi\right)\cdot p_{j,k}(\xi)=p_{j,k-1}(\xi^2)\cdot \alpha(\xi) & \hbox{} \\
    \alpha\left(e^{2\pi i\frac{j}{2^k}}\xi\right)\cdot q_{j,k}(\xi)+\beta\left(e^{2\pi i\frac{j}{2^k}}\xi\right)=
p_{j,k-1}(\xi^2)\cdot \beta(\xi)+q_{j,k-1}(\xi^2)& \hbox{}
  \end{array}
\right.
\]
hence we get exactly the description of the group elements from
Equation \ref{eq:Group Recursive}.

\noindent One can then use the recursive formula to describe each
group element $\gamma_{\frac{j}{2^{k}}}(\xi,z)$. Let $\omega=e^{2\pi i\frac{j}{2^k}}$. Assume that $j$ is odd. Otherwise we would need to look at a smaller $k$. Notice that the first integer $m$ for which $\omega^{2^{m}}=-1$ is $m=k-1$. We compute:
\begin{eqnarray}\label{eq: p}
p_{j,k}(\xi)&=& \frac{\prod\limits_{s=0}^{k-1}\alpha(\xi^{2^s})}{\prod\limits_{s=0}^{k-1}\alpha((\omega\xi)^{2^s})}\\
q_{j,k}(\xi)&=&\frac{p_{j,k-1}(\xi^2)\beta(\xi)+q_{j,k-1}(\xi^2)-\beta\left(\omega\xi\right)}
{\alpha\left(\omega\xi\right)}\nonumber
\end{eqnarray}
\noindent We choose the standard trivialization as in Proposition \ref{Prop:triv1} such that
\[
\lim\limits_{a\rightarrow 0}\frac{\alpha(\xi)}{a}=\frac{\gamma(\xi)}{2\gamma^{2}(\xi^{2})}.
\]
Then $\beta(\xi)=-\alpha(\xi)$ and the relation for $q_{j,k}$ is
\begin{eqnarray}\label{eq: q}
q_{j,k}(\xi)=\frac{-\sum\limits_{s=0}^{k-1}\prod\limits_{t=s}^{k-1}\alpha\left(\xi^{2^t}\right)+\sum\limits_{s=0}^{k-1}\prod\limits_{t=s}^{k-1}\alpha\left((\omega\xi)^{2^t}\right)}{\prod\limits_{s=0}^{k-1}\alpha\left((\omega\xi)^{2^s}\right)}.
\end{eqnarray}
\noindent Define for simplicity $\Pi_{s}(\xi) = \prod\limits_{t=s}^{k-1}\alpha\left(\xi^{2^t}\right)$ for $0\leq s\leq k-1$. So $\Pi_{s}(\omega\xi) = \prod\limits_{t=s}^{k-1}\alpha\left((\omega\xi)^{2^t}\right)$ and in particular $\Pi_{k-1}(\xi)=\alpha\left(\xi^{2^{k-1}}\right)$ and $\Pi_{k-1}(\omega\xi)=\alpha\left(-\xi^{2^{k-1}}\right)$. The formulas for $p_{j,k}$ and $q_{j,k}$ simplify to
\begin{eqnarray}
p_{j,k}(\xi)&=& \frac{\Pi_{0}(\xi)}{\Pi_{0}(\omega\xi)} \label{Eq:p}\\
q_{j,k}(\xi)&=&\frac{\Pi_{k-1}(\omega\xi)-\Pi_{k-1}(\xi)}{\Pi_{0}(\omega\xi)}+\sum\limits_{s=0}^{k-2}\frac{\Pi_{s}(\omega\xi)-\Pi_{s}(\xi)}{\Pi_{0}(\omega\xi)}. 
\label{Eq:q}
\end{eqnarray}

Let $\delta := \inf_{\xi \in \s^{1}}|\gamma(\xi)|$. Suppose $p(x)=x^{2}+c$ is hyperbolic with connected Julia set $J_{p}$ and assume $\gamma:\s^{1}\rightarrow J_{p}$ is the Carath\'eodory loop of $p$. The critical point $x=0$ is in the interior of the filled-in Julia set $K_{p}$ so
\[
    0<\delta \leq |\gamma(\xi)|\leq 2.
\]
Moreover $p(\gamma(\xi))=\gamma(\xi^{2})$ and
$p(\gamma(-\xi))=\gamma(\xi^{2})$. This gives $\gamma(\xi)^{2} =
\gamma(-\xi)^{2}$ and $\gamma(-\xi)=-\gamma(\xi)$. Note that
$\gamma(\xi)$ is not equal to $\gamma(-\xi)$ since otherwise the
external rays corresponding to $\xi$ and $-\xi$ land at the same
point $\gamma(\xi)\in J_{p}$ and they are mapped under $p$ to the
same external ray landing at $\gamma(\xi^{2})$. This is possible
only if $\gamma(\xi)=0$, the critical point, which is a
contradiction since $0\in \mathring{K_{p}}$.

\begin{lemma}\label{Lemma:1} There exists $\delta'>0$ such that for all $a$ with $|a|<\delta'$, we have
\[
    \bigg{|}\frac{\alpha(\xi)}{a}\bigg{|}<\frac{2}{\delta^{2}},
\]
for all $\xi\in \s^{1}$.
\end{lemma}
\proof We have that
\[
\lim\limits_{a\rightarrow 0}\bigg{|}\frac{\alpha(\xi)}{a}\bigg{|}=\bigg{|}\frac{\gamma(\xi)}{2\gamma^{2}(\xi^{2})}\bigg{|}.
\]
Fix $\epsilon=1/\delta^{2}>0$. Then there exists $\delta'>0$ such that for all $|a|<\delta'$,
\[
\bigg{|}\bigg{|}\frac{\alpha(\xi)}{a}\bigg{|}-\frac{|\gamma(\xi)|}{2|\gamma(\xi^{2})|^{2}}\bigg{|}<\epsilon,
\]
and in particular
\[\bigg{|}\frac{\alpha(\xi)}{a}\bigg{|}< \epsilon + \frac{|\gamma(\xi)|}{2|\gamma(\xi^{2})|^{2}}\leq\frac{1}{\delta^{2}} + \frac{2}{2\delta^{2}}= \frac{2}{\delta^{2}}.
\]
\qed

\begin{lemma}\label{Lemma:2} There exists $\delta''>0$ such that for all $a$ with $|a|<\delta''$, we have
\[
\bigg{|} \frac{\alpha(\xi)}{a}-\frac{\alpha(-\xi)}{a}\bigg{|}>\frac{\delta}{8},
\]
for all $\xi\in\s^{1}$.
\end{lemma}
\proof We have
\[
\lim\limits_{a\rightarrow 0}\bigg{|}\frac{\alpha(\xi)}{a}-\frac{\alpha(-\xi)}{a}\bigg{|}=\bigg{|}\frac{\gamma(\xi)}{2\gamma^{2}(\xi^{2})}-\frac{\gamma(-\xi)}{2\gamma^{2}(\xi^{2})}\bigg{|} = \bigg{|}\frac{\gamma(\xi)}{\gamma^{2}(\xi^{2})}\bigg{|},
\]
since $\gamma(-\xi)=-\gamma(\xi)$. Fix $\epsilon=\delta/8>0$.  There exists $\delta''>0$ such that for all $|a|<\delta''$,
\[
\bigg{|}\bigg{|}\frac{\alpha(\xi)}{a}-\frac{\alpha(-\xi)}{a}\bigg{|}-\frac{|\gamma(\xi)|}{|\gamma(\xi^{2})|^{2}}\bigg{|}<\epsilon,
\]
and in particular
\[
\bigg{|}\frac{\alpha(\xi)}{a}-\frac{\alpha(-\xi)}{a}\bigg{|}> -\epsilon + \frac{|\gamma(\xi)|}{|\gamma(\xi^{2})|^{2}}\geq -\frac{\delta}{8} + \frac{\delta}{4}= \frac{\delta}{8}.
\]
\qed

\begin{prop}[\textbf{Growth estimate}]\label{prop:Growth} Suppose $j$ is odd. There exists $a_{0}>0$ such that for all $0<|a|<a_{0}$ there exists a positive integer $k_{0}$ such that for all $k\geq k_{0}$
\[
|p_{j,k}(\xi)z+q_{j,k}(\xi)|> \frac{\delta^{3}}{32}\left(\frac{\delta^{2}}{2|a|}\right)^{k-1} - |z|.
\]
The integer $k_{0}$ depends only on $a_{0}$ and $z$.
\end{prop}
\proof From Equations \ref{Eq:q}  and \ref{Eq:p} we get
\begin{eqnarray*}
|p_{j,k}(\xi)z+q_{j,k}(\xi)|  &=& \bigg{|} \frac{\Pi_{k-1}(\omega\xi)-\Pi_{k-1}(\xi)}{\Pi_{0}(\omega\xi)}+\sum\limits_{s=0}^{k-2}\frac{\Pi_{s}(\omega\xi)-\Pi_{s}(\xi)}{\Pi_{0}(\omega\xi)} + \frac{\Pi_{0}(\xi)}{\Pi_{0}(\omega\xi)}z\bigg{|} \nonumber\\
&\geq& \bigg{|} \frac{\Pi_{k-1}(\omega\xi)-\Pi_{k-1}(\xi)}{\Pi_{0}(\omega\xi)} \bigg{|}  -  \bigg{|} \sum\limits_{s=0}^{k-2}\frac{\Pi_{s}(\omega\xi)-\Pi_{s}(\xi)}{\Pi_{0}(\omega\xi)} + \frac{\Pi_{0}(\xi)}{\Pi_{0}(\omega\xi)}z\bigg{|} \nonumber\\
&\geq&  \frac{\big{|}\Pi_{k-1}(\omega\xi)-\Pi_{k-1}(\xi)\big{|}}{\big{|}\Pi_{0}(\omega\xi)\big{|}}  -   \sum\limits_{s=0}^{k-2} \frac{\big{|}\Pi_{s}(\xi)\big{|}+\big{|}\Pi_{s}(\omega\xi) \big{|}}{\big{|}\Pi_{0}(\omega\xi)\big{|}} - \frac{\big{|}\Pi_{0}(\xi)\big{|}}{\big{|}\Pi_{0}(\omega\xi)\big{|}}|z|
\end{eqnarray*}

\noindent {\bf The leading term.}
\begin{eqnarray*}
\frac{\big{|}\Pi_{k-1}(\omega\xi)-\Pi_{k-1}(\xi)\big{|}}{\big{|}\Pi_{0}(\omega\xi)\big{|}} &=& \frac{\bigg{|}\alpha\left(-\xi^{2^{k-1}}\right)-\alpha\left(\xi^{2^{k-1}}\right)\bigg{|}}{\prod\limits_{t=0}^{k-1}\bigg{|}\alpha\left((\omega\xi)^{2^t}\right) \bigg{|}} = \frac{1}{|a|^{k-1}}\cdot \frac{\bigg{|}\frac{\alpha\left(-\xi^{2^{k-1}}\right)}{a}-\frac{\alpha\left(\xi^{2^{k-1}}\right)}{a}\bigg{|}}{\prod\limits_{t=0}^{k-1}\bigg{|}\frac{\alpha\left((\omega\xi)^{2^t}\right)}{a} \bigg{|}} \\
&\geq&\frac{\frac{\delta}{8}}{|a|^{k-1}}\cdot \frac{1}{\prod\limits_{t=0}^{k-1}\bigg{|}\frac{\alpha\left((\omega\xi)^{2^t}\right)}{a} \bigg{|}} \
\end{eqnarray*}

\noindent {\bf The $s$-term.}
\begin{eqnarray*}
\frac{|\Pi_{s}(\xi)|}{|\Pi_{0}(\omega \xi)|}  &=& \frac{\bigg{|}\prod\limits_{t=s}^{k-1}\alpha\left(\xi^{2^t}\right)\bigg{|}}{\bigg{|}\prod\limits_{t=0}^{k-1}\alpha\left((\omega\xi)^{2^t}\right) \bigg{|}} = \frac{\prod\limits_{t=s}^{k-1}\bigg{|}\frac{\alpha\left(\xi^{2^t}\right)}{a}\bigg{|}|a|^{k-s}}{\prod\limits_{t=0}^{k-1}\bigg{|}\frac{\alpha\left((\omega\xi)^{2^t}\right)}{a}\bigg{|}|a|^{k}} =\frac{1}{|a|^{s}}\cdot \frac{\prod\limits_{t=s}^{k-1}\bigg{|}\frac{\alpha\left(\xi^{2^t}\right)}{a}\bigg{|}}{\prod\limits_{t=0}^{k-1}\bigg{|}\frac{\alpha\left((\omega\xi)^{2^t}\right)}{a}\bigg{|}}\\
&\leq& \frac{1}{|a|^{s}}\cdot \frac{\left(\frac{2}{\delta^{2}}\right)^{k-s}}{\prod\limits_{t=0}^{k-1}\bigg{|}\frac{\alpha\left((\omega\xi)^{2^t}\right)}{a}\bigg{|}} =  \left(\frac{2}{\delta^{2}}\right)^{k}\cdot \left(\frac{\delta^{2}}{2|a|}\right)^{s}\cdot \frac{1}{\prod\limits_{t=0}^{k-1}\bigg{|}\frac{\alpha\left((\omega\xi)^{2^t}\right)}{a}\bigg{|}}
\end{eqnarray*}

\noindent Similarly we can show that
\begin{eqnarray*}
\frac{|\Pi_{s}(\omega\xi)|}{|\Pi_{0}(\omega \xi)|}
&\leq&   \left(\frac{2}{\delta^{2}}\right)^{k}\cdot \left(\frac{\delta^{2}}{2|a|}\right)^{s}\cdot \frac{1}{\prod\limits_{t=0}^{k-1}\bigg{|}\frac{\alpha\left((\omega\xi)^{2^t}\right)}{a}\bigg{|}}
\end{eqnarray*}

\noindent Putting together all inequalities we get that
\begin{eqnarray*}
|p_{j,k}(\xi)z+q_{j,k}(\xi)|
&\geq&  \frac{1}{\prod\limits_{t=0}^{k-1}\bigg{|}\frac{\alpha\left((\omega\xi)^{2^t}\right)}{a}\bigg{|}}\left(  \frac{\frac{\delta}{8}}{|a|^{k-1}} - 2\sum\limits_{s=0}^{k-2} \left(\frac{2}{\delta^{2}}\right)^{k} \left(\frac{\delta^{2}}{2|a|}\right)^{s}  - \left(\frac{2}{\delta^{2}}\right)^{k}|z| \right)\!.
\end{eqnarray*}

We can compute explicitly the sum in the middle and get
\begin{eqnarray*}
\sum\limits_{s=0}^{k-2} \left(\frac{2}{\delta^{2}}\right)^{k} \left(\frac{\delta^{2}}{2|a|}\right)^{s}
&=&\left(\frac{2}{\delta^{2}}\right)^{k}\sum\limits_{s=0}^{k-2} \left(\frac{\delta^{2}}{2|a|}\right)^{s} =\left(\frac{2}{\delta^{2}}\right)^{k}\frac{\left(\frac{\delta^{2}}{2|a|}\right)^{k-1}-1}{\frac{\delta^{2}}{2|a|}-1}\\
&=& \frac{2}{\delta^{2}}\cdot\frac{2|a|}{\delta^{2}-2|a|}\cdot\frac{1}{|a|^{k-1}} -\frac{2|a|}{\delta^{2}-2|a|} \left(\frac{2}{\delta^{2}}\right)^{k}\\
&=&C_{1}(a)\frac{1}{|a|^{k-1}} -C_{2}(a)\left(\frac{2}{\delta^{2}}\right)^{k}\!,
\end{eqnarray*}
where $C_{1}(a) := \frac{2}{\delta^{2}}\frac{2|a|}{\delta^{2}-2|a|}$ and $C_{2}(a) := \frac{2|a|}{\delta^{2}-2|a|}$ are constants that depend on $a$.
We get
\begin{eqnarray*}
|p_{j,k}(\xi)z+q_{j,k}(\xi)|
&\geq&  \frac{1}{\prod\limits_{t=0}^{k-1}\bigg{|}\frac{\alpha\left((\omega\xi)^{2^t}\right)}{a}\bigg{|}}\left( \left(\frac{\delta}{8} - 2C_{1}(a)\right)\frac{1}{|a|^{k-1}} +\left( 2C_{2}(a)-|z| \right)\left(\frac{2}{\delta^{2}}\right)^{k} \right)\\
&>&  \frac{1}{\prod\limits_{t=0}^{k-1}\bigg{|}\frac{\alpha\left((\omega\xi)^{2^t}\right)}{a}\bigg{|}}\left( \left(\frac{\delta}{8} - 2C_{1}(a)\right)\frac{1}{|a|^{k-1}} -|z|\left(\frac{2}{\delta^{2}}\right)^{k} \right)\!,
\end{eqnarray*}
since $C_{2}(a)$ is positive. Note that $C_{1}(a)$ can be made arbitrary small. In particular, if $ |a|<\frac{\delta^{2}}{2}\cdot\frac{\delta^{3}}{\delta^{3}+64}$ then $\frac{\delta}{8}-2C_{1}(a)>\frac{\delta}{16}$.
To see this, notice that
\[
\frac{2|a|}{\delta^{2}}<\frac{\delta^{3}}{\delta^{3}+64}\ \ \ \mbox{and so} \ \ \ \frac{\delta^{2}}{2|a|}-1>\frac{64+\delta^{3}}{\delta^{3}}-1 = \frac{64}{\delta^{3}}.
\]
Then
\[
C_{1}(a) = \frac{2}{\delta^{2}}\cdot\frac{1}{\frac{\delta^{2}}{2|a|}-1}<\frac{2}{\delta^{2}}\frac{\delta^{3}}{64} = \frac{\delta}{32}.
\]
and $ \frac{\delta}{8}-2C_{1}(a)>\frac{\delta}{8} -\frac{2\delta}{32} = \frac{\delta}{16}$. Under this assumption we have shown that
\[
|p_{j,k}(\xi)z+q_{j,k}(\xi)| > \frac{1}{\prod\limits_{t=0}^{k-1}\bigg{|}\frac{\alpha\left((\omega\xi)^{2^t}\right)}{a}\bigg{|}}\left( \frac{\delta}{16}\frac{1}{|a|^{k-1}} -|z|\left(\frac{2}{\delta^{2}}\right)^{k} \right)\!.
\]
If $ |a|<\frac{\delta^{2}}{2}\cdot\frac{\delta^{3}}{\delta^{3}+64}$ then $ \frac{\delta^{2}}{2|a|}>\frac{\delta^{3}+64}{\delta^{3}}>1$. Thus there exists an integer $k_{0}$ such that for all $k\geq k_{0}$ we have
\begin{eqnarray}\label{Eq:k0}
\frac{\delta}{16}\frac{1}{|a|^{k-1}} -|z|\left(\frac{2}{\delta^{2}}\right)^{k}>0\ \Leftrightarrow\ \left(\frac{\delta^{2}}{2|a|}\right)^{k-1}>\frac{32}{\delta^{3}}|z|.
\end{eqnarray}
In view of Lemma \ref{Lemma:1} we get
\begin{eqnarray*}
\prod\limits_{t=0}^{k-1}\bigg{|}\frac{\alpha((\omega\xi)^{2^t})}{a}\bigg{|}\leq \left(\frac{\delta^{2}}{2}\right)^{-k}\!.
\end{eqnarray*}
Hence, for $k\geq k_{0}$ we have
\begin{eqnarray*}
|p_{j,k}(\xi)z+q_{j,k}(\xi)|
&>& \frac{1}{\prod\limits_{t=0}^{k-1}\bigg{|}\frac{\alpha\left((\omega\xi)^{2^t}\right)}{a}\bigg{|}}\left( \frac{\delta}{16}\frac{1}{|a|^{k-1}} -|z|\left(\frac{2}{\delta^{2}}\right)^{k} \right)\\
&>& \left(\frac{\delta^{2}}{2}\right)^{k} \left( \frac{\delta}{16}\frac{1}{|a|^{k-1}} -|z|\left(\frac{2}{\delta^{2}}\right)^{k} \right) = \frac{\delta^{3}}{32}\left(\frac{\delta^{2}}{2|a|}\right)^{k-1} -|z|
\end{eqnarray*}
The constant $a_{0}>0$ can be taken to be
\[ a_{0} = \min\left(\delta',\delta'', \frac{\delta^{2}}{2}\frac{\delta^{3}}{\delta^{3}+64}\right)\!,
\]
where $\delta'>0$ and $\delta''>0$ are the same constants from Lemmas \ref{Lemma:1} and  \ref{Lemma:2}.
\qed

\begin{prop}\label{prop:discont}
Let $(\xi_{0},z_{0})\in \s^{1}\times \C$. There exists a neighborhood $U\subset \s^{1}\times \C$ of $(\xi_{0},z_{0})$ such that
\[
\gamma_{\frac{j}{2^{k}}}(U)\cap U = \emptyset
\]
for all elements $\gamma_{\frac{j}{2^{k}}}\in \Gamma_{p,a}$, with $\gamma_{\frac{j}{2^{k}}}\neq id$.
\end{prop}
\proof Fix $a$ with $|a|<a_{0}$ as in the proof of Proposition \ref{prop:Growth}.
Consider a neighborhood $U_{0}\subset \C$ of $z_{0}$ defined by $|z-z_{0}|<\frac{\delta^{3}}{32}$. Then for all $z\in U_{0}$ we have $|z|<|z_{0}|+\frac{\delta^{3}}{32}$. There exists a smallest positive integer $k_{0}$ which is large enough so that for all $k\geq k_{0}$ the following inequality holds
\begin{equation*}
\left(\frac{\delta^{2}}{2|a|}\right)^{k-1}>\frac{64}{\delta^{3}}\left(|z_{0}|+\frac{\delta^{3}}{32}\right) > \frac{32}{\delta^{3}}|z|,
\end{equation*}
for all $z\in U_{0}$.
Hence Condition \ref{Eq:k0} in the proof of Proposition \ref{prop:Growth} is satisfied for all $z\in U_{0}$. For $j$ odd and $k\geq k_{0}$ we get
\begin{eqnarray*}
|p_{j,k}(\xi)z+q_{j,k}(\xi)-z_{0}| &\geq& |p_{j,k}(\xi)z+q_{j,k}(\xi)| - |z_{0}| >  \frac{\delta^{3}}{32}\left(\frac{\delta^{2}}{2|a|}\right)^{k-1} -|z| - |z_{0}|\\
&>&\frac{\delta^{3}}{32}\left(\frac{\delta^{2}}{2|a|}\right)^{k-1} - 2|z_{0}|-\frac{\delta^{3}}{32}\\
&>&\frac{\delta^{3}}{32}\frac{64}{\delta^{3}}\left(|z_{0}|+\frac{\delta^{3}}{32}\right)-2|z_{0}|-\frac{\delta^{3}}{32}  =\frac{\delta^{3}}{32}.
\end{eqnarray*}
It follows that $|p_{j,k}(\xi)z+q_{j,k}(\xi)-z_{0}| > \frac{\delta^{3}}{32}$ for all $z\in U_{0}$ and for all $k\geq k_{0}$ and  $j$ odd.  If $j$ is even then the situation is similar by considering a smaller $k$.  Suppose $k_{0}>2$.  If $1<k<k_{0}$ then we look at the first component of $\gamma_{\frac{j}{2^{k}}}(\xi,z)$ where
\[
\gamma_{\frac{j}{2^{k}}}\pvec{\xi}{z}=\pvec{e^{2\pi
i\frac{j}{2^k}}\xi}{p_{j,k}(\xi)z+q_{j,k}(\xi)}\!.
\]
Set $V_{0} =\{\xi \in \s^{1}: |\arg(\xi)-\arg(\xi_{0})|<1/2^{k_{0}}\}$ and let $U=V_{0}\times U_{0}$ be a neighborhood of $(\xi_{0},z_{0})$. Then
\[
\gamma_{\frac{j}{2^{k}}}(U)\cap U = \emptyset \ \ \ \mbox{for all}\ \ j,k\  \mbox{with} \ \frac{j}{2^{k}}\neq 1.
\]
To summarize, when $k$ is large (i.e. $k\geq k_{0}$)  the second component of $\gamma_{\frac{j}{2^{k}}}(\xi,z)$ exits $U$ and when $k$ is small, the first component of $\gamma_{\frac{j}{2^{k}}}(\xi,z)$ exits $U$. \qed

\section{Main results}\label{sec: Results}
The growth estimates described in Section \ref{sec:grgamma} provide a powerful tool for analyzing the properties of the extension of the group of deck transforms $\Gamma_{p,a}$ to the boundary of the covering manifold. We first recall some basic properties of group actions. 
\begin{defn}\label{def:discont}
Let $X$ be a locally compact metric space. A discrete group $G$ acts
properly discontinuously on $X$ if for every $x\in X$ there exists a
neighborhood $U\subset X$ of $x$ such that $g(U)\cap U = \emptyset$
for all group elements $g\in G$, $g\neq id$.
\end{defn}
\begin{defn}\label{def:free}
Let $X$ be a locally compact metric space. A discrete group $G$ acts
freely on $X$ if for every $x\in X$ $g(x)\neq x$ for all group
elements $g\in G$, $g\neq id$.
\end{defn}

 We are now able to prove the first theorem about the action of the group $\Gamma_{p,a}$ on
 $\s^{1}\times \C$, the boundary of $(\C-\D)\times \C$.
\begin{thm}\label{thm:discont}
Let $p$ be a hyperbolic quadratic polynomial with connected Julia
set. There exists $a_{0}>0$ such that for all $a$ with $0<|a|<a_{0}$
the group $\Gamma_{p,a}$ acts freely and properly discontinuously on
$\s^{1}\times \C$.
\end{thm}
\proof 
It is easy to see that the action of
$\Gamma_{p,a}$ is free on $\s^{1}\times \C$. Take any element
$\gamma_{\frac{j}{2^{k}}}\in \Gamma_{p,a},
\gamma_{\frac{j}{2^{k}}}\neq id$. Then
\[
\gamma_{\frac{j}{2^{k}}}(\xi,z)=(\omega \xi,
p_{j,k}(\xi)z+q_{j,k}(\xi)),\ \mbox{where}\ \omega=e^{2\pi i
\frac{j}{2^k}}\neq 1.
\]
So the only group element that fixes the first component is the
identity.  The fact that the action of the group on $\s^1\times\C$ is
properly discontinuous is the hard part of the theorem and it
follows from Proposition \ref{prop:discont}, which  in turn uses the growth estimates proved in Proposition \ref{prop:Growth} in an essential way.
\qed

\begin{cor}\label{cor: groupQuotient}
$\s^1\times \C \big{/}\Gamma_{p,a}$  and $(\C-\D)\times \C
\big{/}\Gamma_{p,a}$ are topological manifolds, with fundamental
group $\Z[1/2]/\Z$.
\end{cor}
\proof By Theorem \ref{thm: HOV-escapingSet},
$(\C-\overline{\D})\times \C$ is a covering space of $U^+$, so the
action of $\Gamma_{p,a}$ on $(\C-\overline{\D})\times \C$ is
properly discontinuous and without fixed points. By Theorem
\ref{thm:discont}, the action of $\Gamma_{p,a}$ on $\s^1\times \C$
is also properly discontinuous and without fixed points. \qed

In general, it would be interesting to study whether the group always acts properly discontinuous on $\s^1\times\C$ as in Theorem \ref{thm:discont} or whether there are examples of \He maps for which the group has limit sets on $\s^1\times\C$.

\begin{thm}\label{thm: main1}
Let $p$ be a hyperbolic quadratic polynomial with connected Julia
set. There exists $a_{0}>0$ such that for all $a$ with
$0<|a|<a_{0}$ the following hold
\begin{itemize}
\item[(a)] There exists a continuous surjective map $\widehat{\pi}$,  holomorphic on the leaves of the foliation of $J^+$,  that makes the following diagram commute
\[
\diag{\s^1\times \C \big{/}\Gamma_{p,a}}{\s^1\times \C
\big{/}\Gamma_{p,a}}{J^+}{J^+}{\widetilde{H}}{\widehat{\pi}}{\widehat{\pi}}{H}
\]
\item[(b)] There exists a continuous surjective map $\widehat{\pi}$, biholomorphic on $(\C-\overline{\D})\times \C
\big{/}\Gamma_{p,a}$ and holomorphic on the leaves of the foliation of $J^+$,  that makes the diagram commute
\[
\diag{(\C- \D)\times \C \big{/}\Gamma_{p,a}}{(\C-\D)\times \C
\big{/}\Gamma_{p,a}}{\overline{U}^+}{\overline{U}^+}{\widetilde{H}}{\widehat{\pi}}{\widehat{\pi}}{H}
\]
\end{itemize}
\end{thm}
\proof
By Equation \ref{group}, the map $\widetilde{H}:(\C-\D)\times\C\rightarrow (\C-\D)\times \C$ 
satisfies the relation
\[
\widetilde{H}\circ\gamma_{\frac{j}{2^{k+1}}}=\gamma_{\frac{j}{2^{k}}}\circ\widetilde{H}\
\mbox{and}\ \gamma_1=\mbox{id}.
\]
In Lemma \ref{lemma: pi}, we constructed the function
$\pi:(\C-\D)\times\C\rightarrow U^+\cup J^+$, with the property that
$\pi\circ \widetilde{H}=H\circ \pi$. Since $\Gamma_{p,a}$ is a group
of deck transforms, we have $\pi\circ \gamma =\pi$ for any
$\gamma\in \Gamma_{p,a}$. So both $\widetilde{H}$ and $\pi$ descend
to the quotient $(\C-\D)\times \C \big{/}\Gamma_{p,a}$. By abuse of
notation, we will still use $\widetilde{H}$ in place of
$\widehat{\widetilde{H}}$. 
 The map $\widehat{\pi}$ is a continuous
surjection, holomorphic on the leaves on the foliation of $J^+$ and
$U^+$, and $\widehat{\pi}:(\C-\overline{\D})\times \C
\big{/}\Gamma_{p,a} \rightarrow U^+$ is injective, by Theorem
\ref{thm: HOV-escapingSet}. \qed

\begin{remarka}
Part (a) of Theorem \ref{thm: main1} can be viewed as a two
dimensional analog of the Carath\'eodory loop from one dimensional
dynamics. The universal object in this case is not the circle
$\s^1$, but rather a 3-dimensional topological manifold, isomorphic to a
quotient of $\s^1\times\C$ by a discrete group action.
\end{remarka}

\begin{cor}\label{cor:Mc} Let $p$ be a quadratic polynomial with an attractive fixed point.
There exists $a_{0}>0$ such that for all
$0<|a|<a_{0}$ the closure of the escaping set $U^+$ of the \He map $H_{p,a}$ satisfies  $\overline{U}^+\simeq \
(\C-\D)\times \C/\Gamma_{p,a}$. The Julia set $J^+$ is a topological
manifold and $J^+\simeq \ \s^{1}\times \C/\Gamma_{p,a}$.
\end{cor}
\proof By Theorem \ref{thm: LyubichRobertson}, the boundary of the
primary component is homeomorphic to $\s^1$. The projection
$\widehat{\pi}:\s^{1}\times \C/\Gamma_{p,a}\rightarrow J^{+}$ from
Theorem \ref{thm: main1} is bijective. \qed

In one-dimensional dynamics, W. Thurston  \cite{Th} has constructed topological models for the Julia sets of quadratic polynomials as quotients of the unit circle. 
Consider a hyperbolic polynomial $p(x)=x^2+c$ with connected Julia
set $J_p$ and let $\gamma:\s^1\rightarrow J_p$ be the Carath\'eodory loop of $p$. 
Thurston defined an equivalence relation on
$\s^1$ using the  Carath\'eodory loop, $\xi_1\sim\xi_2$ whenever $\gamma(\xi_1)=\gamma(\xi_2)$, and showed that $\s^1/_{\sim}$ is homeomorphic to the Julia set $J_p$.

Similarly, we will introduce an equivalence relation on
$\s^1\times\C$, and respectively on $\s^1\times\C
\big{/}{\Gamma_{p,a}}$ for $a$ small. When there is no confusion, we
will denote the group $\Gamma_{p,a}$ by $\Gamma$
and the orbit of a point $(\xi,z)$ under the group $\Gamma_{p,a}$ by 
\[
\bigO_{\Gamma}\left((\xi,z)\right)=\left\{\gamma_{\frac{j}{2^k}}(\xi,z)\ :\ k\geq0,\ 1\leq
j\leq2^k\right\}.
\]

\begin{defn}[\textbf{Equivalence of points}] Let $\xi_{1},\xi_{2}\in \s^{1}$ and $z\in \C$.
We will say that
\[
(\xi_{1},z)\sim_{p}(\xi_{2},z) \mbox{ if }
\gamma(\xi_{1})=\gamma(\xi_{2}).
\]
\end{defn}
The following elementary proposition will be useful.
\begin{prop}\label{prop: gamma-pol} Let $\xi_{1},\xi_{2}\in \s^{1}$ such that
$\gamma(\xi_{1})=\gamma(\xi_{2})$.
Let $\omega_{1}=e^{2 \pi i
\frac{j}{2^{k}}}$ be a dyadic root of unity, where $j$ is odd. There exists $m$ odd such that 
if we set $\omega_{2}=e^{2 \pi i \frac{m}{2^{k}}}$, then
$\gamma(\omega_{1}\xi_{1})=\gamma(\omega_{2}\xi_{2})$.
\end{prop}
\proof By induction on $k$. We use the fact that
$p(\gamma(\xi))=\gamma(\xi^{2})$, for any $\xi\in \s^{1}$. 
\qed

\begin{prop}\label{prop: equivalent-elements} Let $\xi,\xi_2\in\s^1$ and
$z\in\C$ such that $(\xi_1,z)\sim_{p}(\xi_2,z)$. Let $k$ be a
non-negative integer and $j$ an odd number with $1\leq j\leq 2^k$. There
exists $m$ odd, $1\leq m\leq 2^k$ such that
\[
\gamma_{\frac{j}{2^k}}(\xi_1,z)\sim_p
\gamma_{\frac{m}{2^k}}(\xi_2,z).
\]
\end{prop}
\proof Let $\omega_1=e^{2 \pi i \frac{j}{2^{k}}}$. By Lemma
\ref{prop: gamma-pol}, there exists $m$ odd such that, if
$\omega_2=e^{2 \pi i \frac{m}{2^{k}}}$ then
$\gamma(\omega_1\xi_1)=\gamma(\omega_2\xi_2)$. We will look at the
group elements
\[
\gamma_{\frac{j}{2^{k}}}\pvec{\xi_1}{z}=\pvec{\omega_1\xi_1}{p_{j,k}(\xi_1)z+q_{j,k}(\xi_1)}
\ \ \mbox{and}\ \
\gamma_{\frac{m}{2^{k}}}\pvec{\xi_2}{z}=\pvec{\omega_2\xi_2}{p_{m,k}(\xi_2)z+q_{m,k}(\xi_2)}\!.
\]
Using Equations \ref{eq: p} and \ref{eq: q}, and Lemma
\ref{prop:Lamin1} we get that
\[
p_{j,k}(\xi_1)z+q_{j,k}(\xi_1)=p_{m,k}(\xi_2)z+q_{m,k}(\xi_2),
\]
so the result follows.
\qed

\begin{cor}\label{cor: orb} Let $(\xi_{1},z),(\xi_{2},z)$ be two points in $\s^{1}\times \C$ such that
$(\xi_{1},z)\sim_{p}(\xi_{2},z)$. Then any point in
$\bigO_{\Gamma}\left((\xi_{1},z)\right)$ is $\sim_p$ equivalent to
some other point in $\bigO_{\Gamma}\left((\xi_{2},z)\right)$.
\end{cor}

\begin{remarka}\label{remark: non-canonical} It is not in general true that if $(\xi_1,z)\sim_p
(\xi_2,z)$ and $\gamma_{\frac{j}{2^k}}$ is an element of the group $\Gamma$, then $\gamma_{\frac{j}{2^k}}(\xi_1,z)\sim_p
\gamma_{\frac{j}{2^k}}(\xi_2,z)$, so one cannot extend canonically
the action of the group $\Gamma$ to the space $\s^1\times\C/_
{\sim_p}$. However, by Corollary \ref{cor: orb}, the orbits of
$\s^1\times \C$ under the group $\Gamma$ preserve the equivalence
relation $\sim_p$.
\end{remarka}

One can extend the notion of $\sim_p$ equivalence
to group orbits as follows:
\begin{defn}[\textbf{Equivalence of orbits}]\label{def: orbitEquivalence}
Let $(\xi,z),(\xi',z')\in \s^{1}\times\C$.
We say that
\[
\bigO_{\Gamma}\left((\xi,z)\right)\sim_{p}
\bigO_{\Gamma}\left((\xi',z')\right)
\]
if there exists $(\xi'',z)\in \bigO_{\Gamma}\left((\xi',z')\right) $
such that $(\xi,z)\sim_{p}(\xi'',z)$.
\end{defn}

\medskip
We need to prove first that $\sim_p$ is well defined. Let $(t,y)$
be another point in $\bigO_{\Gamma}\left((\xi,z)\right)$,
\[
(t,y)=\gamma_{\frac{j}{2^k}}(\xi,z), \mbox{ for some } k\geq 0
\mbox{ and } 1\leq j\leq 2^k, \ j \mbox{ odd}.
\]
We will show that one can find $(t',y)$ in $\bigO_{\Gamma}\left((\xi',z')\right)$ such
that $(t,y)\sim_p (t',y)$. We know that there exists $(\xi'',z)$ in
$\bigO_{\Gamma}\left((\xi',z')\right) $ such that
$(\xi,z)\sim_{p}(\xi'',z)$. By Lemma \ref{prop:
equivalent-elements}, there exists $m$ odd, $1\leq m\leq 2^k$, such
that $\gamma_{\frac{j}{2^k}}(\xi,z)\sim_p
\gamma_{\frac{m}{2^k}}(\xi'',z)$. Then
$(t',y)=\gamma_{\frac{m}{2^k}}(\xi'',z)$ is the element that we
want.

Let us show that $\sim_p$ is an equivalence relation. The
fact that  $\sim_p$ is reflexive and transitive is obvious, so we
only show symmetry. 

The symmetry property follows almost directly
from Lemma \ref{prop: equivalent-elements}. Suppose that
$\bigO_{\Gamma}\left((\xi,z)\right)\sim_{p}
\bigO_{\Gamma}\left((\xi',z')\right)$.
There exists $(\xi'',z)$ in $\bigO_{\Gamma}\left((\xi',z')\right) $,
$(\xi'',z)=\gamma_{\frac{j}{2^k}}(\xi',z')$ such that
$(\xi,z)\sim_{p}(\xi'',z)$. Then
$(\xi',z')=\gamma_{\frac{2^k-j}{2^k}}(\xi'',z)$, so by Lemma
\ref{prop: equivalent-elements}, there exists  $m$ odd, $1\leq m\leq
2^k$, such that $\gamma_{\frac{2^k-j}{2^k}}(\xi'',z)\sim_p
\gamma_{\frac{m}{2^k}}(\xi,z)$. Therefore
$
\bigO_{\Gamma}\left((\xi',z')\right)\sim_{p}
\bigO_{\Gamma}\left((\xi,z)\right),
$
which shows that $\sim_p$ is symmetric.

\begin{thm}\label{thm: JpC}Let $p$ be a hyperbolic quadratic polynomial with connected Julia set $J_{p}$. Let $\mathcal{M}_{p,a}:=\s^1\times \C \big{/}\Gamma_{p,a}$.
There exists $a_0>0$ such that for all $a$ with $0<|a|<a_{0}$ there exists a conjugacy
$\widehat{\pi}$ which makes the following diagram commutative
\begin{equation*} \diag{\mathcal{M}_{p,a}/_{\sim_p}}{\mathcal{M}_{p,a}/_
{\sim_p}}{J^+}{J^+}{\widehat{H}}{\widehat{\pi}}{\widehat{\pi}}{H}
\end{equation*}

\noindent where  $\sim_p$ is the equivalence relation from Definition \ref{def:
orbitEquivalence}.
\end{thm}
\proof The function $\widehat{\pi}:\s^1\times \C \big{/}\Gamma_{p,a}
\rightarrow J^+$ from  Theorem \ref{thm: main1} is a continuous
surjection, biholomorphic on the leaves of the lamination of $J^+$.
If the polynomial from which we perturb is $p(x)=x^{2}+c$ and the
parameter $c$ is not chosen from the interior of the main cardioid of the Mandelbrot
set, then the projection $\widehat{\pi}$ is not yet injective.

The projection function $\pi : \s^1\times \C \rightarrow J^+$ was
first constructed in Lemma \ref{lemma: pi}. We defined $\pi(\xi,z)$
as $\pi_{\xi}(z)$, where $\pi_{\xi}$ is the unique biholomorphic map
from $\C$ into $\mathcal{F}_{\xi}$ with the property that
$\pi_{\xi}(0)=c_{0}(\xi)$ and $\pi_{\xi}(1)=c_{-1}(\xi)$. We 
show that if $\pi(\xi_1,z_1)=\pi(\xi_2,z_2)$ then
$\bigO_{\Gamma}(\xi_1,z_1)\sim_p \bigO_{\Gamma}(\xi_2,z_2)$.

Assume therefore that $\pi(\xi_1,z_1)=\pi(\xi_2,z_2)$ for some
points $(\xi_1,z_1)$ and $(\xi_2,z_2)$ from $\s^1\times\C$. Then
$\mathcal{F}_{\xi_1}$ and $\mathcal{F}_{\xi_2}$ represent the same
leaf of the lamination of $J^+$ and the functions
$\pi_{\xi_1}:\C\rightarrow \mathcal{F}_{\xi_1}$ and
$\pi_{\xi_2}:\C\rightarrow \mathcal{F}_{\xi_2}$ are potentially
different parametrizations of the same leaf. The primary component
of the critical locus intersects $\mathcal{F}_{\xi_2}$ at the points
 $c_0(\omega \xi_2)$, where $\omega$ is a dyadic root of unity.
There exists $\omega=e^{2 \pi i \frac{j}{2^k}}$ some dyadic root of
unity such that
\begin{eqnarray}\label{eq: critical}
c_0(\xi_1)=c_0(\omega \xi_2).
\end{eqnarray}
However, the identifications of the primary component of the
critical locus are completely described in Theorem \ref{thm:
LyubichRobertson}, namely we have
\begin{eqnarray}\label{eq: identifications}
c_0(\zeta_1)=c_0(\zeta_2) \mbox{ for } \zeta_1,\zeta_2\in\s^1
\Leftrightarrow \gamma(\zeta_1)=\gamma(\zeta_2).
\end{eqnarray}
From Relations \ref{eq: critical} and \ref{eq: identifications} it
follows that $\gamma(\xi_1)=\gamma(\omega\xi_2)$. The
Carath\'{e}odory loop $\gamma$ verifies the conjugacy relation
$\gamma(\xi^2)=p(\gamma(\xi))$ so we must also have $
\gamma(\xi_1^2)=\gamma(\omega^2 \xi_2^2). $ Therefore, by Relation
\ref{eq: identifications}, the following equality holds true
\[
c_{-1}(\xi_1)=H^{-1}(c_0(\xi_1^2))=H^{-1}(c_0(\omega^2\xi_2^2))=c_{-1}(\omega\xi_2).
\]
Then $\pi_{\xi_1}$ and $\pi_{\omega \xi_2}$ represent the same
parametrization of the leaf,  so we get
\[
\pi(\xi_1,z)=\pi(\omega\xi_2,z), \mbox{ for any } z\in\C.
\]
In particular this gives $\pi(\xi_1,z_1)=\pi(\omega \xi_2,z_1)$. The
projection $\pi$ when restricted to $\omega \xi_2 \times \C
\rightarrow \mathcal{F}_{\xi_2}$ is injective. However, since
$\pi(\xi_1,z_1)=\pi(\xi_2,z_2)$ by hypothesis, we already know which
point from $\omega \xi_2 \times \C$ projects to $\pi(\xi_1,z_1)$.
This is $\gamma_{\frac{j}{2^k}}(\xi_2,z_2)$, where
$\gamma_{\frac{j}{2^k}}\in \Gamma$ is the deck transform
corresponding to the dyadic root of unity $\omega=e^{2 \pi i \frac{j}{2^k}}$. In
conclusion $\gamma_{\frac{j}{2^k}}(\xi_2,z_2)$ and $(\omega \xi_2,
z_1)$ must coincide. By Definition \ref{def: orbitEquivalence} it
follows that $\bigO_{\Gamma}(\xi_1,z_1)\sim_p
\bigO_{\Gamma}(\xi_2,z_2)$. 
\qed

\begin{remarka}
It would be possible to identify $\mathcal{M}_{p,a}/_{\sim_p}$ with
a quotient of $J_p\times \C$ by an equivalence relation induced by
the group orbits of $\Gamma$ on $\s^1\times\C$, using Remark
\ref{remark: non-canonical}.
\end{remarka}

Theorem \ref{thm: JpC} was proven in the context of
perturbations of a hyperbolic polynomial with connected Julia
set. The main ingredients were Theorems \ref{thm: BS-lamination JU}
and \ref{thm: LyubichRobertson}, out of which the first one
is non-perturbative. The description of the critical locus from
Theorem \ref{thm: LyubichRobertson} may also hold throughout
the entire hyperbolic component of the \He connectedness locus that
contains perturbations of a hyperbolic quadratic polynomial with
connected Julia set. The result of Theorem \ref{thm: JpC} could also be extended to this region.

\section{Extension to semi-parabolic H\'enon maps}\label{sec:remarks}

Another extension concerns \He maps with a semi-parabolic
fixed point, which come from perturbations of a polynomial with a parabolic fixed point.

\begin{defn}\label{def:semi-parabolic}
A fixed point $(x,y)$ of $H$ is called semi-parabolic if the
derivative $DH_{(x,y)}$ has two eigenvalues $|\mu|<1$ and
$\lambda=e^{2\pi i p/q}$.
\end{defn}

The set of parameters $(c,a)\in\C^{2}$
for which the \He map $H_{c,a}$ has a fixed point with one
eigenvalue a root of unity $\ds \lambda$, is a curve of equation
\begin{equation*}\label{eq:P-lambda}
\mathcal{P}_{\lambda}:=\left\{(c,a)\in\C^{2}\ | \
c=(1+a)\left(\frac{\lambda}{2}+\frac{a}{2\lambda}\right)-
\left(\frac{\lambda}{2}+\frac{a}{2\lambda}\right)^{2}\right\}.
\end{equation*}

In \cite{RT1} we studied \He maps with a semi-parabolic fixed point and small Jacobian (see also \cite{RT2} for a discussion on a larger class of hyperbolic \He maps). 
Denote by $V=\D_R\times\D_R$  the polydisk from the Hubbard filtration of $\C^2$ depicted in Figure \ref{fig:filtration}.

\begin{thm}[\cite{RT1}]\label{thm:Parabolics}
Let $p(x)=x^{2}+c_{0}$
be a polynomial with a parabolic fixed point of multiplier $\lambda=e^{2 \pi i p/q}$. There exists $a_{0}>0$ such that for
all $(c,a)\in \mathcal{P}_{\lambda}$ with $0\leq|a|<a_{0}$ the \He map
$H_{c,a}$ has connected Julia set $J$ and there exists a
homeomorphism $\Phi: J_p\times \D_R \rightarrow J^{+}\cap V$ such that
the diagram
\[
\diag{J_p\times \D_R}{J^{+}\cap V}{J_p\times \D_R}{J^{+}\cap V}
{\Phi}{\sigma}{H_{c,a}}{\Phi}
\]
commutes. The function $\sigma$ is given by
$ \sigma(\xi,z)=\left(p(\xi),\xi+\frac{a}{2\xi}z\right)$.
\end{thm}

It follows from Theorem \ref{thm:Parabolics} that for \He maps $H_{c,a}$ which are small
perturbations of the parabolic polynomial $p$ inside the parabola $\mathcal{P}_{\lambda}$, the set $J^+$ inside the polydisk $V$ is a
trivial fiber bundle over $J_p$ with fibers biholomorphic to $\D_R$. 
Notice also that  the vertical disks
$\zeta\times \D_R$, $\zeta\in J_p$ that appear in the description of
$J^+\cap V$ correspond to local stable manifolds of points from
the Julia set $J$.
The proof of
Theorem \ref{thm:Parabolics} from \cite{RT1} also implies that the foliation of
$U^+$ and the lamination of $J^+$ fit together continuously. Therefore Theorem \ref{thm: BS-lamination JU}, which was known for hyperbolic maps, also holds true for this class of semi-parabolic \He maps.

The same arguments as in \cite{LR} can be used to prove parts (a),
(b), (d) and (e) of Theorem \ref{thm: LyubichRobertson}, when $H_{c,a}$ is a
small perturbation inside $\mathcal{P}_{\lambda}$ of a quadratic
polynomial with a parabolic fixed point. The reasoning is similar, because the critical point of a quadratic polynomial $p$ with a parabolic or an attracting fixed point or cycle belongs to the interior of the filled-in Julia set $K_p$. So there exists a primary
component of the critical locus $\mathcal{C}_0$ inside $U^+\cap
U^-$ asymptotic to the $x$-axis and there exists a biholomorphic extension of the function $\varphi^+$ from $\mathcal{C}_0\cap V^+$ to 
$\varphi^+:\mathcal{C}_0\rightarrow \C-\D$. The polydisk $V$ can be used
as a trapping region for $\mathcal{C}_0$ when $a$ is small. Since
$J$ is connected from Theorem \ref{thm:Parabolics}, the boundary $\partial \mathcal{C}_0$ of the primary component
 is contained in $J^+$. The proof of part
(c) from Theorem \ref{thm: LyubichRobertson} follows from Theorem \ref{thm:Parabolics}.

\begin{lemma}\label{lemma:CritPar}
Let $p$ be a quadratic polynomial with a parabolic fixed point of multiplier
$\lambda=e^{2 \pi i p/q}$. There exists $a_0>0$ such that for
all $(c,a)\in \mathcal{P}_{\lambda}$ with $0\leq |a|<a_0$ the boundary
$\partial \mathcal{C}_0$ of the primary component $\mathcal{C}_0$ of
the critical locus for the \He map $H_{c,a}$ is homeomorphic to the
Julia set $J_p$ of the parabolic polynomial $p$.
\end{lemma}

Therefore Theorem \ref{thm: JpC} can be generalized to the semi-parabolic setting as follows

\begin{thm} Let $p$
be a quadratic polynomial with a parabolic fixed point of multiplier
$\lambda=e^{2 \pi i p/q}$.  For all $(c,a)\in \mathcal{P}_{\lambda}$ with $0<|a|<a_0$, the group $\Gamma_{c,a}$ acts properly discontinuously  and without fixed points  on $\s^1\times
\C$. Let $\mathcal{M}_{c,a}=\s^1\times
\C/\Gamma_{c,a}$. There exists
a homeomorphism $\widehat{\pi}$ which makes the following diagram commute
\begin{equation*}
\diag{\mathcal{M}_{c,a}/_{\sim_{p}}}{\mathcal{M}_{c,a}/_{\sim_{p}}}{J^+}{J^+}
{\widehat{H}_{c,a}}{\widehat{\pi}}{\widehat{\pi}}{H_{c,a}}
\end{equation*}
\end{thm}
\proof In view of Lemma \ref{lemma:CritPar} and Theorem \ref{thm:Parabolics}, the proof is the same as that of Theorem \ref{thm: JpC}.
\qed 

Notice also that the cocyle $\alpha$ studied in Sections \ref{Affine}
 and \ref{Degeneracy} is a full invariant of the family $\mathcal{P}_{\lambda}$ when the Jacobian is small, because this family is parametrized by the eigenvalue $\mu$ with $|\mu|<1$ of the semi-parabolic fixed point and $\alpha(1)$ equals $\mu$ by Condition \ref{eq:restriction}.

\end{document}